\documentclass[11pt]{article}

\usepackage{makeidx}
\usepackage{geometry}                
\usepackage{graphicx}
\usepackage{amssymb}
\usepackage{epstopdf}
\usepackage{lscape}
\usepackage{pdflscape}
\usepackage{verbatim}

\usepackage[gen]{eurosym}

\newcommand{\pf}{\vs \noindent {\it Proof:} \quad}

\newcommand{\C}{\mbox{\msbm{C}}}

\newcommand{\Q}{\mbox{\msbm{Q}}}
\newcommand{\R}{\mbox{\msbm{R}}}
\newcommand{\Z}{\mbox{\msbm{Z}}}

\renewcommand{\Re} {{\rm Re\,}}
\renewcommand{\Im} {{\rm Im\,}}

\newcommand{\be}{\begin{equation}}
\newcommand{\ee}{\end{equation}}
\newcommand{\ba}{\begin{array}}
\newcommand{\ea}{\end{array}}
\newcommand{\bea}{\begin{eqnarray}}
\newcommand{\eea}{\end{eqnarray}}
\newcommand{\bean}{\begin{eqnarray*}}
\newcommand{\eean}{\end{eqnarray*}}
\newcommand{\1}{{\bf 1}}
\newcommand{\vs}{\vspace{.25in}}

\font\msbm=msbm10

\newtheorem{thmN}{Theorem}
\newtheorem{propN}{Proposition}

\newtheorem{lemN}{Lemma}

\newcommand{\bc}{\begin{center}}
\newcommand{\ec}{\end{center}}
\newcommand{\aster}{\begin{center} *** \end{center}}

\begin{document}

\title{
 Tessellations and Positional Representation}
\vspace{22pt}
\author{
{\sc Howard L. Resnikoff}\footnote{Resnikoff Innovations LLC; howard@resnikoff.com.}
}
\date{20140420}          
\maketitle

\begin{abstract}
The main goal of this paper is to define a 1-1 correspondence between between substitution tilings constructed by inflation and the arithmetic of positional representation in the underlying real vector space. 

It introduces a generalization of inflationary tessellations to equivalence classes of tiles. Two tiles belong to the same class if they share a defined geometric property, such as equivalence under a group of isometries, having the same measure, or having the same `decoration'. Some properties of ordinary tessellations for which the equivalence relation is congruence with respect to the full group of isometries are already determined by the weaker  relation of equivalence with respect to equal measure. In particular, the multiplier for an inflationary tiling (such as a Penrose aperiodic tiling) is an algebraic number.

Equivalence of tiles under measure facilitates the  investigation of  properties of tilings that are independent of dimension, and provides a method for transferring tilings from one dimension to another.

Three well-known aperiodic tilings  illustrate aspects of the correspondence:   a tiling of  Ammann, a  Penrose tiling, and the monotiling of Taylor and Socolar-Taylor.

\end{abstract}

\noindent {\sc Keywords:} {\it Ammann tilings, aperiodic tilings, Fibonacci numbers, golden number,  inflation, measure-preserving maps, multi-radix, Penrose tilings, positional notation, positional representation, non integral radix, remainder sets, silver numbers, substitution tilings, Taylor monotile, tessellations, 2-dimensional positional representations.}


\section{Introduction}

{\it Mathematics in Civilization} \cite{MIC} argued that there are only two problems in mathematics: {\it improving the ability to calculate} and {\it understanding the geometrical nature of space}. As knowledge increases, these fundamental problems are reformulated in a more sophisticated way, and investigated anew. The goal of this paper is to highlight the connection between positional representation for numbers and  geometrical  tilings of the plane.\footnote{Tilings are also called {\it tessellations}, from the Latin {\it tessera} for the small pieces of stone, glass or ceramic tile used in mosaics. `Tessera' is derived from the Greek for `four' referring to the four sides of rectangular mosaic stones. We use `tessellation' and `tiling' interchangeably.} 

Although positional representation for real numbers is one of the most ancient and greatest intellectual inventions  and certainly amongst the most important in numerical practice,  the concept seems not to have had much influence on the internal development of mathematics. Other good, old, mathematical ideas, such as the unique factorization of integers into primes and the theorem of Pythagoras, have led to innumerable deep and important abstract generalizations. Not so positional representation. Positional representation is still the occasional subject of mathematical papers, most often in  the recreational category. 

Tessellations  have an even  more ancient history but they have generally been considered decorative rather than profound.  Although vast numbers of papers and online examples have been devoted to tilings, the mathematical theory is recent and unstructured -- largely a collection of interesting and sometimes beautiful patterns.

Positional representation and tessellations have  traditionally been considered independent domains of recreational mathematics. Connections between the two have  not received much notice.  Potential relationships could be of particular interest for  tilings because they might provide natural constructions as well as an arithmetic representation for the geometric relationships, and an alternative but familiar language for  constructing and  talking about tilings. In the opposite direction, tilings provide insight into new forms of positional representation -- particularly those that depend on more than one radix.\footnote{The term `radix' refers to the base of a system of positional representation.}

The paper explores this interplay.  The main result is a 1-1 correspondence between a class of positional representations and a class of tilings that includes those constructed by the process of inflation used by Roger Penrose\,\cite{Penrose1} in 1974.

\subsection{A very brief history of positional representation}

It is crucial to distinguish between a {\it positional notation} -- a notation that employs a finite inventory of symbols to identify an arbitrary real number -- and a {\it positional representation}, which is a positional notation whose symbols and structure have meanings -- interpretations --  that are linked to the structure of arithmetic so that the notation can be used for calculation.

This distinction may be worth elaborating.  Suppose that a 1-1 correspondence between the field of real numbers $\R$ and a set $S$ without any structure is given. Suppose further that a notation for the elements of $S$ that employs sequences of symbols drawn from a finite inventory is used to set up a correspondence with the elements of $S$. The set of such sequences can be thought of as a  `positional' method for labeling the elements of $S$;  call it a {\em positional notation} for $S$. These sequences contain no hint of arithmetical properties. The 1-1 correspondence  could be used to map the real numbers onto the sequences of the notation, and then to transfer the field properties of $\R$ to $S$ but there is no assurance that this would result in a practical or efficient method for {\em calculating}, that is, for performing the field operations of addition and multiplication and other operations derived from them 

A {\em positional representation} for $\R$ (or for $\C$) is a special kind of correspondence between real numbers  and sequences of symbols drawn from a finite inventory  such that the symbols and the sequences have meanings that explicitly relate them to the numbers and their properties in a way that facilitates the expression of properties of $\R$ (or $\C$) that are not necessarily properties of other sets.  

For instance, the set of potential sentences of English can be coded as sequences of letters from a finite alphabet. This set has the same power as $\R$. Although the letters of the alphabet (together with the interword space symbol) and sequences of letters are also endowed with a linear order (used, for instance, to organize dictionaries), this positional notation for sentences does not imply an {\em a priori}  `arithmetic' of sentences.

The first positional notation that also was an efficient positional representation was introduced by Akkadian mathematicians more than four thousand years ago.\footnote{See \cite{MIC} for the details  of early systems of numeration, and \cite{Knuth2}, section 4.1, for a brief but excellent overview of the history.} It originally was limited to positive integers but was easily extended to positive numbers smaller than 1. In its earliest realization it lacked a symbol for the zero, which was unreliably denoted by a gap between neighboring digits. The radix was 60 -- probably the largest integer ever systematically and extensively employed as a base for hand calculation. Although 60 has the advantage of many divisors and leads to short expressions for the practical quantities that were of interest in early times, the addition and multiplication tables are too large to be memorized, or even used, by anyone other than a specialist. Remnants of radix 60 representation are found in our notation for angles and time.

Radix 10, often referred to as Arabic notation and less often but more correctly as  Hindu-Arabic notation, is said to have been invented  between the 1st and 4th centuries  by Indian mathematicians.  It was adopted by Arab mathematicians many centuries later and made its way to western Europe during the Middle Ages.  Leonardo Pisano (Leonardo of Pisa), generally known as Fibonacci, brought radix 10 calculation into mainstream european thought in his book {\it Liber abaci} --  {\it The Book of Calculation} \cite{Fibonacci, Sigler} --  first published in 1202. Today few remember Fibonacci's role in the transmission of Hindu-Arabic radix 10 representation of numbers to Europe, but many have heard of the sequence of numbers  his book introduced as the solution to a homework  problem -- the  `Fibonacci numbers'.  This sequence, which begins $1, 1, 2, 3, 5, 8, \dots$ (the $n$-th term is the sum of the previous two), plays a role in  many unexpected places, from the idealized reproduction of rabbits, which was the subject of the exercise,  to the  growth of petals on flowers and seeds on pine cones, and not least of all,  in many of the examples in this paper. The analytical formula for Fibonacci's sequence,   said to have been discovered by the 17-th century French  mathematician Abraham de Moivre, is expressed in terms of $\phi = (1+\sqrt{5})/2 \sim 1.61 $ -- the `golden number', defined by Euclid as the ``extreme and mean ratio", that is, the number  $x$ satisfying $x = 1+1/x$. The larger  solution of $x^2 = x + 1$  is $\phi$.  The golden number was  believed to express aesthetically pleasing proportions for rectangles and  for that reason it is often seen embodied in commercial designs and logos. 

No history of positional representation, no matter how brief, should omit Donald Knuth's introduction of 2-dimensional radices \cite{Knuth1, Knuth2}, several variants of which appear below. These ideas were followed up by many others, in particular Gilbert, who studied arithmetic in complex bases in a series of interesting papers \cite{Gilbert1, Gilbert 1982, Gilbert2}. 

The theory of wavelets, a concept  introduced about thirty years ago, is closely related to positional representation although  the connection has not been emphasized. `Wavelets' are collections of compactly supported orthonormal  functions that, in general, overlap and are  bases for a broad variety of function spaces \cite{HLR+ROW 2}. One can think of compactly supported wavelets as a kind of positional representation for functions. Their supports progressively decrease in size as they `home in' on the neighborhood of an arbitrary point on the line.  We shall not examine wavelets here, but the reader should be aware of the intimate and unexplored connection of that circle of ideas with positional representations and tessellations.\footnote{Cp. \cite{Lawton-Resnikoff 1991, Strichartz 1993}.} 

There have not been many fundamental applications of positional representation in pure mathematics but there are a few. 

Positional notation was first used to prove significant theorems by the inventor of set theory,  Georg Cantor,  who recognized that the properties of positional {\it notation} alone -- no need for the additional arithmetical implications of positional representation -- were sufficient to prove that the set of real numbers is not countable\,\cite{Cantor diag method}. Cantor's `diagonal method' has become a foundation stone in the education of mathematicians. He did not neglect positional representation: Cantor's construction of a nowhere decreasing continuous function that increases from 0 to 1 but  is constant except on a set of measure zero (`Cantor's function', cp.  \cite{Cantor function})  used the essential device of passing from positional representation with radix 3 to positional representation with radix 2. 

It has long been known that positional representation provides a way to map subsets of $\R^m$ to subsets of $\R^n$. The general method is made clear from the simplest example: mapping the unit interval onto the unit square. Having fixed the radix, say $\rho=2$, from the representation  of $u \in \R$ as $u = \sum_{k \geq 1} u_k 2^{-k}$, construct the pair $(x,y)$ as
$$
x = \sum_{k \geq 1} u_{2k-1}2^{-k}, \quad  y= \sum_{k \geq 1} u_{2k}2^{-k}
$$
After suitable normalization,  $u \rightarrow (x,y)$, and similar maps constructed from positional representations, is measure-preserving.  Norbert Wiener made essential use of this in his generalized harmonic analysis.\footnote{\cite{Wiener NPRT}, esp. p.81. This property plays a  role in constructing measures on spaces of functions in the theory of brownian motion.}

\subsection{A very brief history of tessellations of the plane}

Tessellations of the plane are no doubt older than the first developments of positional representations of numbers, and developed examples can already be found in the ancient Fertile Crescent. An example from the Sumerian city Uruk IV, circa -3100, now in the Pergamon Museum in Berlin, shows triangular and diamond periodic mosaic patterns.  

Tessellations were taken up from a mathematical viewpoint  in 1619 by Kepler, who wrote in his {\it Harmonices Mundi} --  {\it Harmony of the World} --  about coverings of the plane by regular polygons. 

Periodic tilings of the plane can be classified by their symmetries into 17 groups, sometimes called ``wallpaper groups" or, more formally, ``plane crystallographic groups." While the mathematical classification was the work of Evgraf Fedorov \cite{Fedorov} in 1891, the classification has been intuitively understood for  millennia by artists and craftsmen around the globe who decorated almost every surface they could find with complex repeating patterns.

Fedorov also recognized that crystals were physical realizations of periodic tiling of  3-dimensional space. This led to his classification of the 230 space groups -- the symmetry groups of crystallographic tessellations -- which is among the earliest and most significant mathematical results in this field.

It was not until 1974 that the mathematician and mathematical physicist Roger Penrose \cite{Penrose1} discovered the  aperiodic tessellations that bear his name. Since then it has become a parlor game for amateur  and  professional mathematicians to find new and interesting examples of aperiodic tessellations, but the subject has not stimulated much work nor found  resonance in  other departments of mathematics.

In the past, it seems to have been generally believed that tessellation of the plane without periodic symmetry is impossible. The concepts underlying symmetry and their appearance in art and nature as well as their applications in mathematics and science were traced in a beautiful book  written in 1952 by Hermann Weyl \cite{Weyl-symmetry}, who was one of the most powerful mathematical minds of his time. That Weyl made no mention of the possibility of aperiodic tilings demonstrates how improbable they were thought to be.  

\section{Tessellations}

This paper explores the relationship between tessellations and positional representation, primarily  in the real vector spaces $\R$  and $\R^2$  equipped with the euclidean metric and the measure $m$ derived from it, although some concepts are formulated for $\R^d$.  Two measurable subsets of $\R^d$ are {\it essentially disjoint} if the measure of their intersection is 0, and  they are {\it essentially identical} if the measure of their intersection is equal to the measure of each of the sets. Thus  $m(A \cup B) = m(A) + m(B)$ if and only if $A$ and $B$ are essentially disjoint. 

The objects of interest are `tiles' and the tessellations made from them. A {\it tile} is a subset of $\R^d$ that has positive measure. Tiles are usually, but need not be, connected. Each tile  belongs to one of a finite collection of equivalence classes.  Tiles are of the same {\it type} if they belong to the same equivalence class. If two tiles are of the same type, each is a {\it copy} of the other. We suppose that every type of tile has infinitely many members, and that there is at least one type.  A {\it tessellation} of $\R^d$, also called a {\it tiling}, is the pair consisting of a finite collection of types of tiles of dimension $d$ and a covering of $\R^d$  by essentially disjoint  tiles each of which belongs to one of the types. An {\it overtiling}\,\footnote{This concept is used only in a footnote on page \pageref{overtiling}.}  is the pair consisting of a finite collection of types of tiles of dimension $d$ and a covering of $\R^d$  by essentially identical  tiles each of which belongs to one of the types. In an overtiling, tiles can be stacked on top of one another.

Here are some examples of equivalence classes of tiles:  Congruence under the group of isometries of $\R^d$ is an equivalence relation.  Congruent tiles  are of the same type relative to this relation.  Tiles that are congruent under some subgroup of the group of isometries can also be said to constitute a type. Congruent tiles have the same measure. Tiles that are connected and simply connected and have the same measure but not necessarily the same shape form an equivalence class.   Tiles might be distinguished by `decorations'. Suppose a measurable subset of $\R^d$ is given and copies of it are colored in a finite number of distinct colors. Say that two tiles are equivalent if they have the same color. There are as many types of tiles as colors used to color them.   Tiles may have more general decorations, that is, have markings on them.   Those that have the same decoration form an equivalence class and constitute a type. The decorations may be used to limit how adjacent tiles may be placed. For instance, if the decorations are curves drawn on the tiles,  an allowed tessellation might be one for which the curves are continuous across tile boundaries.

Suppose that the dimension is $d$ and  that there are $N$ types of tiles. Let $\{ R_j: 1 \leq j \leq N \}$ be a set of representatives of the types of tiles. A tessellation is said to be {\it inflationary} if there is a real number $\rho>1$ and an orthogonal linear transformation $O$ such that each magnified representative $ \rho  \,O(R_i) $  is the union of essentially disjoint tiles.\footnote{One could generalize this definition to expansive matrices but for our purposes that would only complicate the details.} This process of magnifying each $R_i $ by the same factor and then tiling it with copies of the $R_j$ is the process called {\it inflation}; the factor $\rho \, O$ is called the {\it multiplier}; sometimes we shall refer to $\rho$ itself as the multiplier. The measure of $ \rho  \,O(R_i)$ is $\rho^d m(R_i )$ so repetition of the inflation process covers increasing volumes of $\R^d$. If the origin of the magnification  lies in the interior of  $R_i$, infinite repetition of the inflationary process results in a tessellation of $\R^d$ that is said to have been {\it constructed by inflation}.  If there is an  $S \subset \R^d$ and a lattice $\Lambda \subset \R^d$ (that is, a discrete additive subgroup of  $\R^d$ of rank $d$) such that $\R^d= \bigcup_{\lambda \in \Lambda} \left( \lambda+ S \right) $ and the translations $ S \rightarrow  \lambda+ S$  are essentially disjoint, then the tiling is said to be {\it periodic}; otherwise it is {\it not periodic}. 

A collection of types of tiles for which at least one tessellation is possible but no tessellation consistent with the constraints is periodic is said to be {\it aperiodic}. In an important paper that connected tessellations to undecidability problems, Wang \cite{Wang 1961} conjectured that aperiodic tilings are impossible. Five years later,  Berger \cite{Berger 1966} showed the existence of aperiodic tilings of the plane by creating a correspondence between tilings and  Turing machines  and applying the undecidability of the halting problem. He constructed a set of 20,426 distinct types of tiles for which an associated tessellation exists and is aperiodic. But it  was Penrose's\,\cite{Penrose1} explicit construction of aperiodic tilings using two types of tiles in 1974 that captured the imagination. In both cases the `types' are defined by the equivalence relation of geometrical congruence up to sets of measure zero.

\aster 

We shall see that inflationary tessellations are intimately related to algebraic numbers, and that aperiodicity is the generic situation. But first, consider tessellations that are both periodic and inflationary.

\begin{thmN} 
Suppose that a tessellation of $\R^d$ is both periodic with period lattice $\Lambda$ and inflationary with multiplier $\rho \, O$. Then 
$$  \rho\,O( \Lambda) \subset \Lambda $$
and  $\rho$ is an algebraic number of degree $d$.  After suitable normalizations, if $d=1$ then $\rho \in \Z$; if $d=2$ then $\rho$ is an imaginary quadratic integer.
\label{thm 1}
\end{thmN}
\vspace{-18pt}
\pf    If $\{ \omega_k : 1\leq k \leq d \}$  is an integral basis for $\Lambda$ there is a matrix $A$ with integer entries such that $ \rho \, O( \omega_j ) = A\, \omega_j$.  Then  $\rho^d = | \det A | $ so $\rho$ is an algebraic number of degree $d$. 

If $d=1$ then $ \rho \, \omega = A \, \omega$ with $A \neq 0$ a rational integer. If $d=2$, then $\R^2$ can be identified with $\C$ in the usual way,\footnote{
 The isomorphism is $ x+i y \leftrightarrow \left( \ba{cc}x & y \\
 -y & x \ea \right)$.} $(\rho \, O)$ with a complex number   temporarily denoted $\rho$, and $\rho \, \omega_j = \sum_j A_{ij} \omega_j$ with $A = 
\left(
\ba{cc}
a & b \\
c & d
\ea \right)$ a matrix of integers. $A$ is invertible because $\Lambda$  is non-degenerate. If the tessellation is inflationary with   a non zero rational integer multiplier, say $n$,  then for any lattice  $\Lambda$ the periods of  $n \Lambda$ are $n \omega_1, n \omega_2$ so $n \Lambda  \subset \Lambda$. 

Are there special lattices for which other multipliers exist? Set $\tau=\omega_1/\omega_2$. Without loss of generality, suppose that $\Im(\tau)>0$; then $\rho= c \tau + d $ and $ \tau = \frac{a \tau + b}{c \tau+d}$. Thus $\tau$ is a quadratic algebraic number. Since the lattice has rank 2, the quadratic is irreducible and the number is imaginary quadratic. 

To complete the proof, we show that $\rho$ satisfies $\rho^2 -(a+d) \rho +\det A =0$ by calculating $\tau$. The discriminant of the quadratic for $\tau$  is $(a-d)^2 + 4bc = (a+d)^2 - 4 \det A$. Thus 
\bean
\rho & =& d +c \tau =d + c \left(  \frac{a-d \pm \sqrt{ (a+d)^2 - 4 \det A }}{2c} \right) \\
&=&   \frac{a+d \pm \sqrt{ (a+d)^2 - 4 \det A }}{2}  
\eean
so $\rho^2 -(a+d) \rho + \det A = 0.  $ 
Since $0 \neq \det A \in \Z $  and $\rho = c \tau +d $, it follows that $\rho$ is an imaginary quadratic integer. $\Box$\\

A consequence of this theorem is that periodic inflationary tessellations of the plane correspond to complex multiplication.\footnote{Including the rational integers.}  Since the the Penrose tilings are inflationary and the multiplier is $ \rho = \frac{1+\sqrt{5}}{2}$ -- not an  imaginary quadratic number -- it follows that no particular Penrose tiling can be periodic, whence Penrose tilings are aperiodic.

The theorem has a partial converse.

\begin{thmN}
If  $d \in \{1,2 \}$ and $\Lambda \subset \R^d$ is a lattice such that $  \rho\,O(  \Lambda ) \subset \Lambda$  for some $\rho \, O \in  \Lambda $ then there is an inflationary tiling with multiplier $\rho \, O$.
\end{thmN}
\vspace{-18pt}
\pf If $d=1$ then $\Lambda \subset \R$ can be normalized so that  $\Lambda =  \Z$. Choose a fundamental domain  $F=\Z/ \Lambda $. For example, the interval $F=[0,1]$ is essentially identical to a fundamental domain. Then the sets in $ \bigcup_{n \in \Z} (n + F )  $ are essentially disjoint and the union is a periodic tessellation of $\R$.   According to the hypothesis, there is a multiplier $\rho \in \Z^*$. Evidently $\rho \, \Z \subset \Z$ and $\rho \, \Z = \bigcup_{k \in \Delta} (k + F)$ where $\Delta$ is a set of representatives of $ \Z / \rho \Z$, for instance $\{ 0, 1 \dots, |\rho|-1 \}$. This is the inflationary decomposition of the tile $F$.

If $d=2$ then identify $\R^2$ with $\C$.  $\Lambda$ has two generators which can be normalized to $\{ 1, \tau \}$ where $\tau \in \C$ has positive imaginary part. According to the hypothesis, $\Lambda$ has a multiplier $\rho$ which  is an integer in an imaginary quadratic number field $\Q(\sqrt{D})$ where $D$ is a negative square-free integer. Let $\Z[\sqrt{D}]$ be the ring of integers in this field and let $\Delta$  be a set of representatives for $ \Z[\sqrt{D}] / \rho \, \Z[\sqrt{D}]$ which we call the set of {\it digits} for the radix $\rho$.  The number of elements in  $\Delta$ is $\rho \overline{\rho} = | \rho |^2$ (The ratio of the area of $\rho F$ to the area of $F$ is the number of congruent copies of $F$ that tile $\rho F$). If $F$ is essentially identical to a fundamental domain $\C/ \Lambda$, then the sets in $\bigcup_{\lambda \in \Lambda} (\lambda + F)$ are essentially disjoint and the union is a periodic tessellation of $\C$. The tessellation is inflationary because $\rho \, F = \bigcup_{\delta \in \Delta} (\delta + F) $.  This is nothing more than saying that each $\lambda \in \Lambda$ can be written as $\lambda = \rho \lambda' + \delta$. $\Box$
\vspace{18pt}

Denote the ring of integers in $\Q(\sqrt{D})$ by $\Lambda = \Z[\sqrt{D}]$. It is a lattice generated over $\Z$ by 1 and 
$$
\omega = \left\{
\ba{lll}
\sqrt{D} ,			& 	D \equiv 2,\, 3  &(\bmod \, 4) \\
\frac{1 +  \sqrt{D} }{2}, & 	D  \equiv 1  & (\bmod \, 4)
 \ea
\right.
$$
Thus $\rho = n_1 + n_2 \omega$ with $n_1, n_2 \in \Z$. 

These tilings are primarily of interest because of their connection with the theory of elliptic functions and number theory.  Another reason is that they  yield positional representations for complex numbers with respect to the radix $\rho$, and thereby establish a connection between a class of tessellations of the plane and positional representations with an imaginary quadratic integer radix.  Indeed, suppose that $\rho \in \Lambda$ and $| \rho |^2 >1$,  and write  the inflation relation  in the form
$$ F = \frac{1}{\rho} \bigcup_{\delta \in \Delta} (\delta  +  F)$$ 
where $\Delta$ is a complete set of representatives of $\Lambda / \rho \Lambda$. Then $F$ is just the remainder set and the  positional representation for an arbitrary remainder is 
\be
z = \sum_{k=1}^{\infty} z_k \rho^{-k}, \quad z_k \in \Delta
\ee

When the cardinality of $\Delta$ -- the number of digits in the positional representation --  is 2, there are  three distinct imaginary quadratic fields and three  possible multipliers up to multiplication by a unit of the field. Listed in order of increasing trace of the complex generator of the  field, they are  $ \rho = i \sqrt{2}$, $\frac{1 + i \sqrt{7}}{2}$, and $1+i$. In each case a convenient choice of digits is $\{ 0, 1 \}$ so these can be considered as generalizations to $\C$ of the conventional binary representation for $\R$.  Pictures of the three remainder sets and their decompositions are shown on pages 168-9 of reference \cite{HLR+ROW 2}.

\aster

Here is a simple example of the correspondence between a tiling and a positional representation with $\rho = 2$. There will be $ | \rho |^2=4$ digits. The simplest tiling of the plane by unit squares -- repetition of one type of square as in a cartesian coordinate grid --   is periodic with lattice generated by $z  \mapsto z+1,\, z \mapsto z+i$. The inflationary multiplier is the rational integer 2.  The tiling has an algebraic realization by the inflation equation
\be
2 R = R \cup  \left(1 + R \right) \cup  \left(1 + i+R \right) \cup  \left( i+ R \right)
\label{squares}
\ee
which is satisfied by $R = \{ z=x+i y \in \C: 0  \leq x, y  \leq 1 \}$. With this solution, eq(\ref{squares}) is an essentially disjoint union. In this case multipliers that map the lattice into itself have the form $m + n i, \, m,n \in \Z$, which are imaginary quadratic (``Gaussian")  integers.

Equation\,(\ref{squares})  is also a realization of  a system of positional representation for the field of complex numbers. It realizes the arithmetic of $\C$ within the framework of a positional representation in the sense that $R$  can be considered a remainder 
and every $z \in \C$ can be written in the form $ z = [ z ] + \{ z \}$
where 
$$\{ z \} =  \sum _{k=1}^{\infty} \frac{z_k}{2^k}  \in R, \quad z_k \in \Delta =\{0,1,1+i,i \}  $$
 and   $[ z ] $ is a Gaussian integer. The elements of $\Delta$ are the digits of the representation. This representation is not unique, just as the representation of real numbers by decimals is not unique. Nevertheless,  each $z$ has a representation relative to  the  remainder set. 
 
 Let us call a number with a  positional representation that has only non negative powers of the radix a {\it positional representation integer}. In the example above, the set of positional representation  integers coincides with the ring of Gaussian integers; hence it is a lattice. But it is not always the case that the integers of a number field and the  positional representation integers of an associated positional representation are the same. This difference can sometimes be exploited to prove aperiodicity of a tessellation.
 
\aster

A checkerboard is also a tessellation of the plane by unit squares but now the squares are of two types, distinguished not by a geometrical property but by a decoration -- their color, say black and white.  This tiling can also be described algebraically. If $B$, resp. $W$, denotes a black, resp. white, square then
\be
\ba{lcl}
2 B &=& B \cup  \left(1 + i + B  \right)\cup  \left(1 + W \right)  \cup  \left( i+ W \right )\\
2 W &=&  B \cup  \left(1 + i + B  \right)\cup  \left(1 + W \right)  \cup  \left( i+ W \right)
\ea
\label{checkerboard}
\ee
In this example each set $B$ and $W$ is geometrically similar to a scaled-up version of itself, and can be decomposed into an essentially disjoint union of copies of tiles.

The colors are not geometrical properties of the tiles, but their role can be replaced by manifestly geometrical properties by deforming the tiles. One way would be to deform the boundaries of the squares: for the black, cut out a triangular notch from two adjacent sides and add semicircular pips to the opposite sides so that the area of the square is conserved; for the white, cut out semicircular notches and adjoin triangular pips so that the tiles can be joined as in a jigsaw puzzle. These modifications force certain relationships in the tessellation.  These relationships are also expressed by eq(\ref{checkerboard}).

Yet another way to think about the checkerboard might be to consider the tiles as 2-faced and to `color', or otherwise distinguish, opposite faces. Reflection in the plane of the checkerboard in an ambient $\R^3$ would represent the mapping from one face to the other.

The literature does not seem to have a general theorem claiming that an arbitrary equivalence relation could be replaced by differences in the shapes of the tiles although something like that must be true for a limited category of equivalence relations. The  modified hexagon  of the Taylor monotile shown in fig.\,\ref{Taylor monotile}  (page \pageref{Taylor monotile}; cp.\,\cite{Socolar+Taylor 2011}) is the most complicated example of this process known to the author; in this case the geometrical alteration results in a disconnected tile.

\aster

The next example is a version of the famous Penrose aperiodic tessellation. Here there are also two types of tiles: isosceles triangles denoted $R_{113}$ and $R_{122}$; the subscripts are the multiples of $\pi/5$ that are their angles. Non-translational isometries of the tiles also appear. Denote complex conjugation by an overline and let $u = \exp ( i \pi/5)$ be  the generator of the  group of rotations of order 10. The inflation factor -- the radix $\rho$ -- is the golden number  $\phi = \frac{1+\sqrt{5}}{2}$. The inflation equations are\footnote{These equations were used to generate the Penrose pinwheel shown in fig.\,\ref{Penrose R122 pinwhl lines} on page  \pageref{Penrose R122 pinwhl lines}.
}
\be
\ba{lll}
\phi \, R_{113} &=& \left( 1 + u^4 \, R_{113} \right) \cup  \left( R_{122}  \right) \\
\phi \, R_{122} &=&   \left( u \,  \overline{R_{113}}  \right)  \cup    \left( 1 + u^3 \, R_{122}  \right)  \cup  \left( 1+ u^5 \,  \overline{R_{122}}   \right) 
\ea
\label{Penrose eqs}
\ee
The tiles that correspond to the remainder sets are isosceles triangles. Their decomposition, implied by the inflation equation eq(\ref{Penrose eqs}), is shown in fig.\,\ref{Penrose 1-digit} on page \pageref{Penrose 1-digit}.

These equations lead to generalized positional representations, at first for points in the remainder sets, and then, by inflation and finally rotation of the wedge-shaped sectors, for arbitrary points of the complex plane. The translations appearing in eq(\ref{Penrose eqs}) are the digits: $\Delta =\{ 0, 1 \}$. There are several  remainder sets that are reflected and rotated as the positional representation advances from  digit to digit.  

\begin{figure}[h]
\begin{center}
\caption{Radix $\phi \sim 1.61$. Penrose  $R_{113}$ 8-digit expansions.}
\vspace{6pt}
\label{pix Penrose R113 8-digit orbit}
\includegraphics[width=4in]{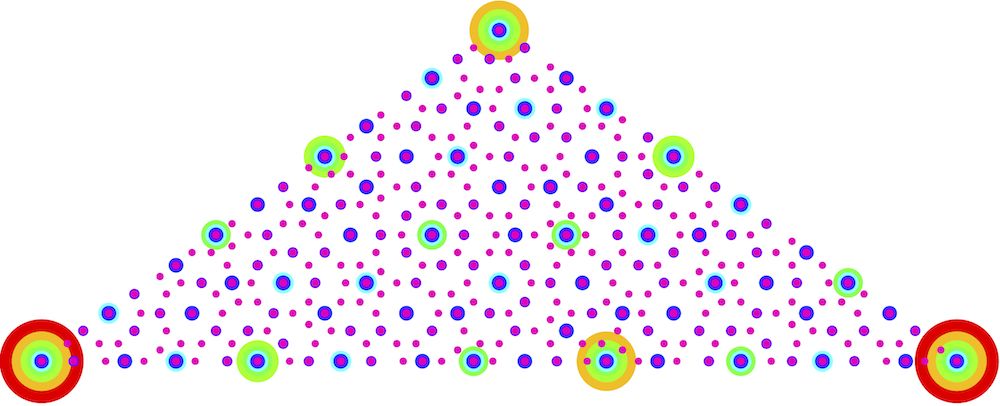}
\end{center}
\end{figure}

\begin{figure}[h]
\begin{center}
\caption{Radix $\phi \sim 1.61$. Penrose  $R_{122}$ 8-digit expansions.}
\vspace{6pt}
\label{pix Penrose R122 8-digit orbit}
\includegraphics[width=4in]{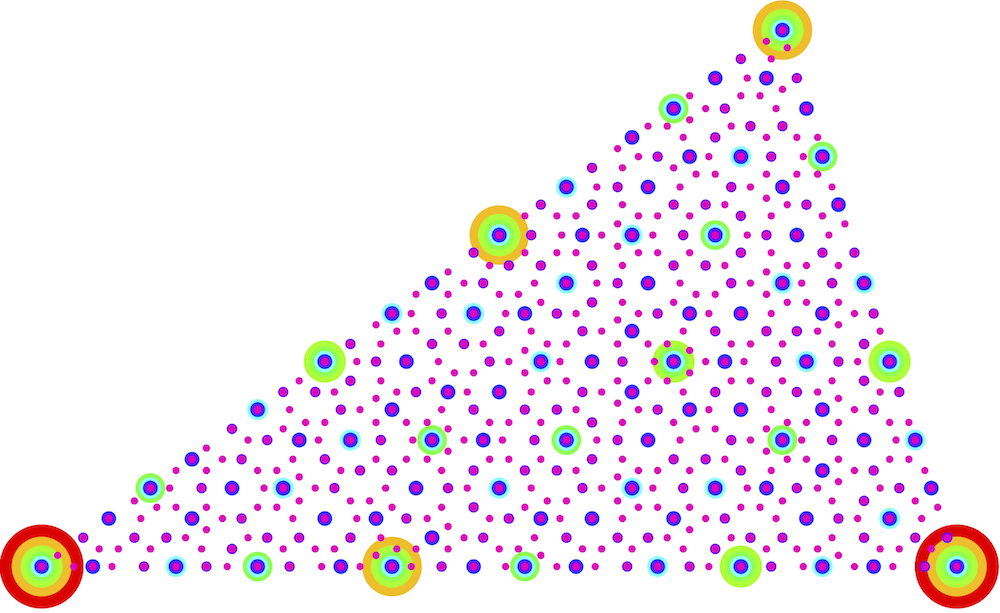}
\end{center}
\end{figure}

The positional representations have the form
\be
z = \sum_{k} \frac{z_k}{\phi^k} u_k, \quad z_k \in \Delta, \, u_k \in \{u^n: 0 \leq n <10 \}
\ee
and each remainder set is the set of all representations of the kind specified by the equations.  Figure\,\ref{pix Penrose R113 8-digit orbit}, resp. fig.\,\ref{pix Penrose R122 8-digit orbit},  displays the expansions in $R_{113}$, resp. $R_{122}$, through 8 digits. The color coding and diameter of the disk that represents a number are arranged so that  the points  corresponding to a given number of digits can easily be seen. The large red disks correspond to $\{ 0/\phi, 1/\phi \}$; the smaller orange disks correspond to 2-digit expansions, etc.

Inflation by $\phi$ extends this to a sector in the plane, and rotation by powers of $u$ extends the sector to the entire plane. Observe that the  triangular remainder sets $R_{113}$ and $R_{122}$  -- which were originally thought of as the  tiles -- are fully determined by the equations, and the partition of each remainder set is essentially disjoint. 
 
The point that corresponds to a finite expansion is a vertex of a deflated copy of a remainder set. In this sense, it labels the deflated remainder set.  In each remainder set the positional representation for a number defines a polygonal path -- call it simply a {\it path} -- from the origin to the point representing the number and labeling the deflated remainder set: just  add the complex numbers -- the vectors --  corresponding to successive digits.  The paths are a microscope that opens up a universe of  detail  in the progressively deflated remainder sets. 

\begin{figure}[h]
\begin{center}
\caption{8-digit path in  $R_{113}$ from 0 to $\frac{1}{\phi } +\frac{u^4}{\phi^2}   +\frac{u^8}{\phi ^3} +\frac{u^{2}}{\phi ^4} + \frac{u^{6}}{\phi ^5}    +\frac{1}{\phi ^6}    +\frac{u^4}{\phi ^7} 
 +\frac{u^8}{\phi^8}$.}
\vspace{6pt}
\label{path R133 6dig}
\includegraphics[width=4in]{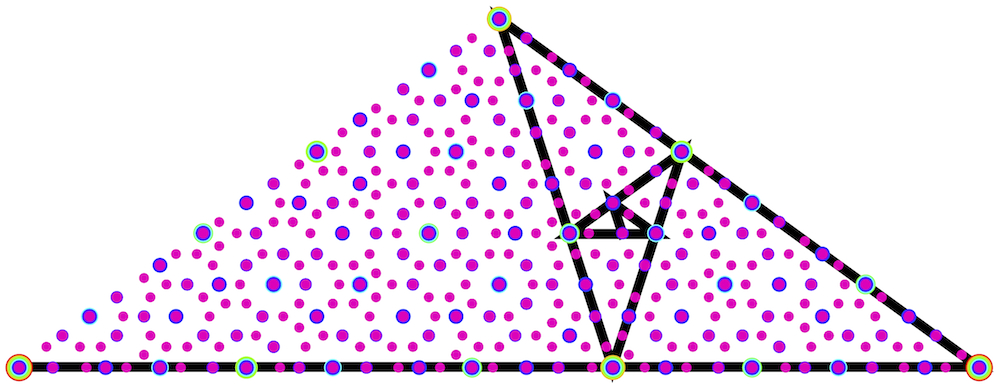}
\end{center}
\end{figure}

For instance, fig.\,\ref{path R133 6dig} (page \pageref{path R133 6dig}) shows an 8-digit path from $0$ to 
$$\frac{1}{\phi } +\frac{u^4}{\phi^2}   +\frac{u^8}{\phi ^3} +\frac{u^{2}}{\phi ^4} + \frac{u^{6}}{\phi ^5}    +\frac{1}{\phi ^6}    +\frac{u^4}{\phi ^7}  +\frac{u^8}{\phi^8} = -\frac{35}{2}+8 \sqrt{5}+\frac{1}{2} i \sqrt{85-38
   \sqrt{5}}$$
    It is easy to trace out the steps digit by digit. The powers of $u$ express the changes of direction along the progressively smaller segments of the path. Note that each turn is by the angle $4\pi/5$.

\aster

It has been said that a remarkable property of Penrose tilings is that every finite pattern of tiles is repeated infinitely often somewhere else in the tiling.  An unremarkable property of positional representations is that every pattern -- every sequence -- of digits in the representation of a number is repeated in the sequence of digits of infinitely many numbers.  For inflationary tilings with positional representations, the discussion above shows that  these facts are the same.

Let us elaborate this observation. Any finite portion $T$ of the tessellation of the plane can be deflated until it becomes a subset of each of the remainder sets. Within the remainder set each tile has a distinguished vertex -- the point in a subtile that corresponds to 0 in the reminder set --   labelled by a finite expansion in the positional representation. Thus the collection of tiles is in a correspondence with a collection of positional representations, and  hence with the corresponding paths. If $T$ is a subset of a deeply nested remainder set $R$, then the initial segments of the paths will coincide until the vertex specifying $R$ is reached, after which they may diverge to the tiles they represent within $R$. Exactly the same is true for the corresponding numbers and their positional representations. So all of these numbers will share an initial sequence of digits $D_{{\rm initial}}$ until some place in the notation, at which point a vertex of the appropriate deflated tile has been specified, and different sequences of digits $\{D_k \}$  will follow the initial common segment $D_{{\rm initial}}$, one for each tile in $T$.  

Let us do some gene splicing. Suppose a new fixed sequence of digits $D_{{\rm new}}$ is inserted between the common segment $D_{{\rm initial}}$ and each $D_k $. This has the geometrical effect of deflating the tiles $T$ further. Looked at through the other end of the microscope, when the tessellation has been re-inflated to the level where the tiles are their original size, they will be located somewhere else, depending on the path for the sequence $D_{{\rm new}}$. 
 \aster
This is the model we will generalize. It suggests that the idea of a positional representation be extended to include a remainder set for each type of tile. Our attention will generally be restricted to ambient spaces $\R$ or $\R^2$ and $\R^2$ will be identified with $\C$. Let $R_i$ be a finite collection of tiles for an inflationary tessellation with radius $\rho$. The inflationary hypothesis is equivalent to
\be
\rho \, R_i = \bigcup_j \left( \delta_{ij} + u_{ij}(R_j) \right), \quad \delta_{ij} \in 
\Delta, \, u_{ij}  \in {\cal O }
\label{PR union}
\ee
where ${\cal O }$ is a finite subgroup of the orthogonal group. In the plane, the action $R \rightarrow u_{ij}(R)$ is either $R \rightarrow u\,R$ or  $R \rightarrow u\,\overline{R}$  where $u$ is a root of unity.
 Iteration of this system of relations shows that for each $i$ and  $z \in R_i$ there are digits $z_k$ such that
\be
z = \sum_{k=1}^{\infty} \frac{z_k}{\rho^k} u_k
\ee
This is a positional representation for $z$ that corresponds to the tessellation.

Thus

\begin{thmN}
If an inflationary tiling corresponds to a positional representation, then every compact part of the tiling  is contained infinitely often in other parts of the tiling. 
\end{thmN}

\vspace{-18pt}
\pf Let us recapitulate what has already been said for the Penrose tiling. After a possible rotation to bring the region into concordance with the remainder sets, every finite collection $T$ of tiles is contained in an inflated version of a remainder set, say $\rho^n R_i$. The tiles in the  collection have a positional integer representation. These leading digits coincide in the deflated set  $\rho^{-n} T \subset R_i$. Select an arbitrary finite sequence of digits and rotations and follow it by the given one. The concatenation is the distinguished vertex of a deflated remainder set after some number of iterations. Inflating that set  produces a region of the tiling in which $T$ reappears.  $\Box$

\aster

Let us restrict our attention to classes of tiles that have the same measure. Consider an inflationary tessellation of $\R^d$ formed from $N$ types of tiles.  An $N$-rowed matrix  $U$, called the {\it partition matrix},  can be constructed from this information. The element $U_{ij}$ is the number of tiles of type $j$ required in the tessellation of the inflated tile $ \rho \, O( R_i )$. This number is a non negative rational integer. The vector $v$ whose $j$-th component is $m(R_j)$ --  the measure of $R_j$ --  is an eigenvector of $U$. Indeed, from the definitions,\footnote{For  rotations $O$ the determinant is 1; for orthogonal transformations with $\det{O}=-1$, package the sign with $\rho$.}
\be
U v = \rho^d v
\label{U.1}
\ee
The  components of the eigenvector lie in the field generated by the eigenvalue $\rho^d $. Conversely, given the eigenvector $v$, the eigenvalue can be expressed as
$$ \rho^d = \frac{v^t U v}{v^t v} $$
which shows that $\rho^d $ lies in the field generated by the measures of the different types of tiles, and $\rho$ lies in that field with  $d$-th roots adjoined.

\begin{thmN}
 A multiplier  $\rho$ for an inflationary tessellation is an algebraic integer of degree at most $N d$.
\end{thmN}
\vspace{-18pt}
\pf $\rho^d$ is a root of the characteristic polynomial of $U$, which has rational integers as coefficients and leading coefficient   1.  Therefore $\rho^d$ is an algebraic integer of degree at most $N$. $\Box$\\

These simple remarks already tell us that not every real number greater than 1 in absolute value can be the multiplier of an inflationary tessellation. This limitation, which interweaves algebraic number theory with geometry,  is part of what makes the class of inflationary tessellations interesting.

If there is only one type of tile then $U$ is a matrix of order 1 whose sole entry is the positive integer $n>0$  that counts how many copies of the tile $R$ are needed in an essentially disjoint tiling of $\rho R$. Hence the eigenvalue equation is $\rho^d = n$ which determines $\rho$ up to a root of unity. This is the case that applies to periodic tilings using a single type of tile, like the  white marble hexagons that were used by interior designers  of an earlier generation  to tile the  bathroom floors. Here the tiles are congruent but the result applies more generally to tiles that have equal area. In particular, if there were a monotile -- a single tile that covers the plane aperiodically in an inflationary tessellation -- this equation would apply. A tiling of the plane by a process similar to inflation that uses just one type of tile was discovered by Joan Taylor \cite{Taylor 2010}; it will be discussed below.\footnote{It is inflationary only  in the limit. Inflation by the multiplier yields an overtiling of $\rho R$. \label{overtiling}}

The partition matrix $U$ contains information that is sometimes sufficient to prove that an associated tessellation is aperiodic. The basic idea goes back to Penrose \cite{Penrose1} but here it appears in the more general setting of equivalence relations not necessarily limited to geometrical congruence.

The multiplier $\rho^d$ is an eigenvalue of the $N \times N$ matrix  $U$.  The entry $U_{ij}$ is the number of copies of tiles of type $j$ required to partition the inflated tile $(\rho \, O) R_i$ of type $i$. In particular, $ \sum_{j=1}^N U_{ij}$ is the total number of tiles needed to tile $\rho R_i$. The partition matrix of the $n$-th iterated inflation is $U^n$; the entries on its $i$-th row  are the number of tiles of each type required to tile $(\rho \, O)^{n} R_i$ and $\sigma^n_i:=\sum_{j=1}^N (U^n)_{ij}$ is the total number of tiles required to tile $(\rho \, O)^{n} R_i$.

\begin{thmN}
 Let $U$ be the partition matrix of a tessellation.  If there is a $j$ such that $\lim_{n \rightarrow \infty} \left(U^n \right)_{ij}/ \sigma^n_i $ is not rational, then the tessellation of $R_i$ with matrix $U$  is aperiodic.
 \end{thmN}
\vspace{-18pt}
\pf Suppose the tessellation were periodic. Then there would exist some set $S$ whose translates form an essentially disjoint cover of $\R^d$. This set can be covered by an integral number of tiles of the various types, so the ratios $ \left(U^n \right)_{ij}/ \sigma^n_i $ are all rational. Take an increasing sequence of essentially disjoint translates: each corresponding ratio for the union is constant, so each limit is also rational.  Thus no tiling with partition matrix $U$ is periodic, so the tiling is aperiodic.  $\Box$ \\

This method of showing aperiodicity will be referred to as the {\it proof by irrationality}. The theorem provides a way to prove that some tilings consisting of tile types that differ only in color are aperiodic. The critical ingredient is a number-theoretic property of the partition matrix. Note that this approach to proving aperiodicity cannot be used if there is only one type of tile.

\aster

Thus far it has been shown that a multiplier for an inflationary tessellation must be an algebraic integer. We have seen an example of an aperiodic tessellation and made a connection between it and a generalization of positional representation where the two types of tiles correspond to two types of remainder sets.

Now we will indicate how this generalization of positional representation arises naturally -- why more than one type of remainder set occurs -- and show that it further constrains the multiplier. 

The simplest way to see this is to try to construct 1-dimensional  inflationary tessellations -- tilings of $\R$ -- with multiplier $1 <\rho \leq 2$. Suppose the tile  $R$ is an interval.  Moreover, consider a  cover of $\rho R$ by two copies of $R$: $\rho R = R \cup (v + R)$ with $v \in \R$. Iteration of this equation yields a positional representation for elements of $R$ as
\be
R = \left\{  x : x= \sum_{k \geq 1} \frac{\delta_k}{\rho^k} v \right\}, \quad \delta_k \in \{0,1 \} 
\label{overlapped intervals}
\ee
and it follows that the endpoints of $R$ are $0$ and $\sum_{k \geq 1}  \rho^{-k} v =  v/(\rho - 1) $. Without loss of generality, normalize the length of $R$ by setting $v=1$. Then $R=[ 0, 1/(\rho -1) ]$.

Observe that $R \cap ( 1 +R) = [1,  1/(\rho -1) ]$ so the two-set cover of $R$  will  be essentially disjoint only when $\rho = 2$ and only then will it produce a tiling. This motivates us to introduce   new types of tiles to make the cover essentially disjoint for other $1 < \rho < 2$.

$\rho R$ can be expressed as the essentially disjoint union $\rho R = [0,1] \cup (1 + R)$, which introduces the new tile $R_1 = [0,1]$. Under inflation $R_1$ is the essentially disjoint union 
$$ \rho R_1 = [0, \rho] = R_1 \cup ( 1+  [0,\rho -1] )$$
This introduces another new tile, $R_2 = [0, \rho -1]$ that inflates as
$$ \rho R_2 = R_1 \cup \left( 1 + [0, \rho (\rho - 1) - 1 ] \right)$$
This process can be repeated indefinitely\footnote{Subject to the obvious inequalities on $\rho$ that insure the number  $x_n$ is not negative.} yielding a sequence of tiles $R_n = [ 0, x_n]$  with 
$$
x_n =\rho \left(\rho \left( \dots  \left( \rho-1\right)  - 1 \right) \dots  - 1\right) -1 = \rho^n -\sum_{k=0}^{n-1} \rho^k
$$
The number of types of tiles can be limited to $n$ by requiring 
$$ \rho^n = \sum_{k=0}^{n-1} \rho^k$$
This constraint makes the space spanned by the $\rho^k$ over $\Q$  $n$-dimensional; forces the multiplier to be an algebraic integer of a special type; and replaces the original overlapped set union for $R$ by a collection of interwoven essentially disjoint set unions for the $R_k$.

This roughly indicates  how positional representations with many remainder sets arise.

\aster

We will fill in the details of this procedure to produce an infinite collection of aperiodic tessellations, starting with aperiodic tessellations of $\R$ and  use them to construct aperiodic tessellations of higher dimensional spaces, concentrating on $\R^2$.

The complexity of both results and presentation increases rapidly as the radix grows. In order to keep the calculations comparatively simple and the discussion informative it will be helpful to concentrate on positional representations that, like the binary system, only use the digits 0 and 1. Constructing these positional representations is the task of section \ref{sec:3}.

\section{Silver numbers}
\label{sec:3}

If $1 < \rho \leq 2$ then $\rho$ can be used as the radix of a positional representation  with digits   $ \{ 0, 1 \}$, which  may be thought of  as a kind of  `generalized binary representation'. Every $x$ in the remainder set $R$  has a positional representation
\be
x  = \sum_{k=1}^{\infty} \frac{x_k}{\rho^k}, \quad x_k \in \{ 0, 1 \}
\label{rho.PR}
\ee 
from which it follows that 
$$R=\{ x: 0 \leq x \leq 1 /(\rho -1)\}$$
The series for $x$  can be conveniently expressed in positional notation as
\be
x = (0 \cdot  x_1 x_2 \dots x_n \dots)_{\rho} 
\label{rho.PN}
\ee
If the radix $\rho$ is fixed,   write
$$x = 0 \cdot x_1 x_2 \dots x_n \dots $$
We make the convention that an infinite repetition of the sequence $x_k\, x_{k+1} \dots \, x_{k+l}$ is abbreviated $\underline{x_{k}\, x_{k+1} \dots \, x_{k+l}}$ and that an infinite sequence of  trailing zeros may be omitted.

\aster
Among these generalized binary representations are a family with particularly interesting  properties.  Let $N \in \Z^+$ and select  a sequence of $N$ `bits' $b_j \in \{0, 1 \} $ with $ b_N=1$. Put $b = (0  \cdot b_{1} \dots  b_N )_2 $, and assume that $2^N b$ is odd and  greater than 1. The integer $2^N b$ satisfies $1 <  2^N b <2^N$.  The {\it silver number} of  index $b$, denoted $s_{b}$, is the largest real root of the polynomial of degree $N$
\be
 P_{b} (x):=x^N\left( 1 -  \sum_{j=1}^{N} b_j x^{-j} \right) 
\label{SN.def}
\ee

\begin{lemN} The silver number $ s_b $ exists and $1< s_b <2$.
\end{lemN}
\vspace{-18pt}
\pf   By definition, $2^N b>1$ is odd. If $x \geq 2$ then 
$$ P_b (x) >  P_b (2) = 2^N \left(1 -\sum_{j=1}^{N} b_{j} 2^{-j} \right) \geq  2^N \left( 1 - \sum_{j=1}^N 2^{-j} \right) >0,$$ while $ P_b (1) = 1-\sum_{j=1}^{N} b_{j} \leq 0 $. Since  $ P_b (x)$ is continuous, it has a  real root, and hence a largest real root, between 1 and 2. Since $ P_b (x) >  P_b (2)$ for $x>2$, it has no larger real root.
$\Box$\\

For $N=2$ there is one silver number:  $s_{3/4}=  \phi =\frac{1+\sqrt{5}}{2}$, the golden number. Note that $\lim_{N \rightarrow \infty}  s_{b} =2$. 

\begin{thmN}
  The silver numbers are algebraic integers of degree greater than 1.
\end{thmN}
\vspace{-18pt}
\pf A silver polynomial is monic with rational integer coefficients. Thus its roots are algebraic integers. The proof reduces to showing that $ P_b (x)$ is irreducible over $\Q$. Suppose otherwise. Then  there exist distinct co-prime integers $p, q$ such that $ P_b (p/q)=0$. Since $b_N=1$, this assumption  implies 
$$p^N - \sum_{j=1}^{N-1} b_{j} p^{N-j} q^{j} = q^N $$
The left side is divisible by $p$  (recall that $N>1$); the right side is not. $\Box$\\

\aster

The polynomial $ P_b $  generates identities that produce distinct finite positional representations for the same number. This comes about as follows.   Let the radix be $\rho= s_b $. Then $\rho$ satisfies \be
1 =(0 \cdot  b_1 \dots b_N)_{\rho}
\label{silver identity}
\ee
This  implies that  there are infinitely many  finite positional representations\footnote{That is, expansions such as eq(\ref{rho.PR}) that have finitely many non-zero  digits.}  that denote the same number. The golden number provides the simplest example. With $b=3/4$, division of eq(\ref{silver identity}) by the appropriate power of $\rho $ shows that the following representations are equal:
$$ 0 \cdot 1 = 0 \cdot  0 1 1 =  0 \cdot  0 1 0 1 1  =  0 \cdot  0 1 0 1 0 1 1  \quad \mbox{etc.}$$

\aster

We take the Frobenius companion matrix of $ P_b(x) $ in the form
\be
U :=  ( u_{ij} ), \qquad u_{ij }=  \left\{
\ba{ll l}
b_j & \mbox{if} & i=1  \\
1  & \mbox{if} & j = i-1 \quad \mbox{and} \quad 1<i \leq N  \\
0 & \mbox{else} &
\ea
\right.
\label{partition matrix}
\ee
that is,

\be
U = 
\left(
\ba{cc c cc}
b_1 & b_2 & b_3 & \dots & b_N \\
1& 0 & 0 & \dots & 0 \\
0& 1 & 0 &  \dots &  0 \\
\vdots & & \ddots && \vdots \\
0 & 0 & \dots & 1 & 0
\ea
\right)
\label{companion matrix Pb}
\ee
where $b_j \in \{ 0,1 \}$.
The characteristic polynomial of $U$ is $ P_b (x)$ so $ s_b $ is the largest real eigenvalue of $U$.

Fix the degree $N$ and the index $b$ and drop them from the notation.  

\begin{lemN}
The silver number $s$ is an eigenvalue of $U$.  If $v$ is an eigenvector of $U$ that belongs to $s$, then  $v$  is proportional to the vector whose  components are $v_k = s^{-k}$. 
\end{lemN}
\vspace{-18pt}
\pf The components of an eigenvector  $v$ for $s$ satisfy $s \,v_k  = v_{k-1} $ for $2 \leq k \leq N$ and $s\, v_1 = \sum_{k=1}^N b_k v_k = 1$. Then $v_k  = c\, s^{-k}$, $c \neq 0$.  Since $v \neq 0$ it can be normalized so that $c=1$, i.e. $\sum_k v_k=1$.   $\Box$ \\

\section{Tessellations of $\R$}

Silver numbers are connected to tessellations by identifying the companion matrix $U$ with the partition matrix for an inflationary tessellation. Properties of the companion matrix thereby become properties of tessellations and, conversely, properties of the partition matrix for a tessellation become properties of the corresponding positional representation.

Every $\rho \in \R$ such that $1< \rho \leq 2$ is the radix of a positional representation whose digits are drawn from $ \{0, 1 \}$. The remainder set is
\be
R = \left\{ x : x = \sum_{k=1}^{\infty} \frac{x_k}{\rho^k} \right\} , \quad x_k \in  \{0, 1 \} 
\label{R def}
\ee
$R$ is an interval one of whose endpoints is 0. The other is $\sum_{k=1}^{\infty} \frac{1}{\rho^k}  =\frac{1}{\rho -1} $, a number that is greater than 1 unless $\rho=2$. Separating the first digit from the sum in eq(\ref{R def}) shows that $R$ satisfies the set theoretic identity
\be
\rho R = R \cup  (1 + R ) 
\label{R decomp}
\ee
As was mentioned above, unless $\rho=2$, the two sets on the right are not essentially disjoint because their intersection contains the interval $[1,1/(\rho-1)]$. This raises the problem of replacing $R$ by a collection of remainder sets such that the decomposition corresponding to eq(\ref{R decomp}) will be a union of essentially disjoint subsets. Here is where the special properties of silver numbers come in.

Suppose that $U$ is the companion matrix of a silver polynomial and  $\rho$ is the largest real eigenvalue of $U$ -- the associated silver number. The eigenvector belonging to to $\rho$ is proportional to $v = (v_k)$ where $v_k = \rho^{-k}$, $1 \leq k \leq N$.  

Equation\,(\ref{R def}) implies that $\rho R= [0,\rho/(\rho-1)]$.  The defining property of the silver number $\rho$ implies 
$$\frac{\rho}{\rho -1}= \frac{\rho}{\rho -1}  \sum_{k=1}^{N} b_k \rho^{-k} $$
which tells us that the length of $\rho R$ can be written as the sum of lengths  of $N$ essentially disjoint intervals.  So introduce sets 
\be
R_k =  \frac{ \rho}{\rho-1} [0,  \rho^{-k} ], \quad 1\leq k \leq N
\label{R_k def}
\ee

The intervals $b_k R_k$ can be arranged in any order as essentially disjoint subsets that cover $R$.  We shall say that the remainder set $b_k R_k$ is {\it essential} if $b_k \neq 0$, and that the order  $R_1, \dots, R_N$ is the {\it natural order}; note that  the lengths of the intervals are  decreasing although only the essential ones contribute to the essentially disjoint covering of $R$. 

Introduce constants 
$$c_k =  \left\{
\ba{lll}
0 & \mbox{if} & k=1 \\
\sum_{j=1}^{k-1} b_j \rho^{-j} & \mbox{if} & 1<k \leq N
\ea
\right.
$$
The remainder sets $R_k$ arranged in the standard order satisfy the system of equations
\be
\ba{lcll }
\rho  \, R_1 &=&   \bigcup_{k=1}^N \left(c_k + b_k R_k \right) \\
\rho  \, R_{k} &=& R_{k-1}, & 2 \leq k \leq N  \\
\label{dim1 rem decomp}
\ea
\ee
 
\begin{thmN}
The decomposition eq(\ref{dim1 rem decomp}) is  an essentially disjoint cover of the original remainder set $R$  and $R = \rho R_1$. Moreover, the decomposition implies the measure relationships expressed by the companion matrix.
\end{thmN}
\vspace{-18pt}
\pf For the first part it is enough to observe that the right hand endpoint of $R_k$ is the left hand endpoint of $R_{k+1}$ (so the intervals are essentially disjoint) and that the sum of their lengths is $\frac{\rho}{\rho-1} \sum_{j=1}^{N} \rho^{-j} =\frac{\rho}{\rho-1}$, which is the length of $\rho R$. 

For the second, the measure is invariant under  translations and reflection   so $m(c_k  +R_k) = m(R_k)$ and the decomposition implies $U v = \rho v$. $\Box$
\vspace{18pt}

This procedure produces an infinite number of `tilings' of $\R$ similar to barcode or a strange piano keyboard. They may be thought uninteresting. However, although they do not look like much to the eye, it turns out that they are, in general, aperiodic. They also provide a foundation for constructing higher dimensional tessellations that are also aperiodic.

\begin{figure}[h]
\begin{center}
\caption{Radix $\phi$:  10-digit decomposition of remainder set $R_2$.}\label{d=1_s(3/4)_R2_n=10_barcode}
\includegraphics[width=5in]{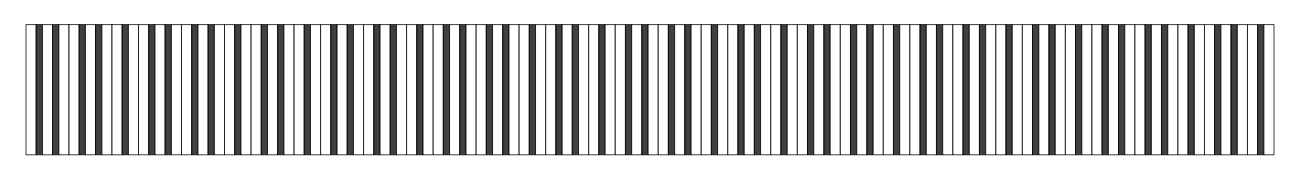}
\end{center}
\end{figure}

For instance, if $\rho =s_{3/4} = \phi$ then the companion matrix is 
given by eq(\ref{Ammann-Penrose U}). Figure \ref{d=1_s(3/4)_R2_n=10_barcode} shows the 10-digit decomposition of $R_2$ constructed as above. It turns out that  $U^2$ is the  matrix for dimension 2 that describes the aperiodic plane  tessellations of Ammann and Penrose.

\aster

These tilings are aperiodic. We shall sketch two methods of proof. The first, following Penrose, shows that the ratio of the number of tiles of the two tile types is irrational. The second proves the non existence of a period lattice.

\begin{thmN}
If $\rho$ is a silver number whose companion matrix $U$ defined by  eq(\ref{companion matrix Pb}), then $U$ is the partition matrix for a  inflationary tessellation of $\R$ and the tessellation is aperiodic.
\label{thm: aperiodic silver}
\end{thmN}
\vspace{-18pt}
\pf The proof can be reduced to a tessellation in standard order because the calculations that will be made are order-independent.

Since $\rho$ is a silver number, $\rho$ is an algebraic irrational and the digits in the associated positional representation are 0 and 1.

The partition matrix elements count the number of tiles belonging to each class that are required to cover each tile representative $\rho R_k$. $\rho R_1$ is partitioned into a total of $\sum_{k=1}^N b_k$ tiles whereas every other $\rho R_k$ is covered by $R_{k-1}$, hence by 1 tile.  After $n$ inflations the number of tiles of type $j$ required to partition $\rho^n R_i$ is the matrix element $\left( U^n \right)_{ij}$.   

$U$ determines a linear recurrence sequence and the matrix entries are the numbers in the sequence. Let us make this explicit.  Let $\1$ denote the vector whose entries are 1. Then the $j$-th entry of $U^n \1$ is the number of tiles required to partition $\rho^n R_j$. Denote this number by $a_{n +j}$. This sequence is generated by the relation $a_{n} = \sum_{j=1}^N  b_j a_{n-j}$, or equivalently, by 
$$(a_{n+1}, \dots  a_{n+1+j}, \dots, a_{n+1+N})^t = U (a_{n}, \dots  a_{n+j}, \dots, a_{n+N})^t$$
It is well-known that the solution has the form
$$ a_n = \sum_{j=1}^N c_j \lambda_j^n $$
where $c_j$ are constants and the $\lambda_j$ are the eigenvalues of $U$, i.e. the roots of $P_b(x)$.  The product of the roots is 1 since   $b_N=1$. Recall that  $\rho>1$. An application of Rouch\'{e}'s theorem shows that each root has absolute value bounded by $\rho$. Hence the eigenvalue $\rho$ dominates and therefore $\lim_{n \rightarrow \infty} a_{n+1}/a_n = \rho$, which is irrational because $P_b(x)$ is irreducible of degree greater than 1.   This implies  the 1-dimensional  tessellation defined by $U$ is aperiodic. $\Box$

\aster

For the golden number, the recurrence defined by the partition matrix $U$  generates the Fibonacci sequence. It follows that if ${\bf e}_j:=(\delta_{jk} )^t$ where $\delta_{jk}$ is the Dirac delta, then $\left( U^d \right)^n {\bf e}_j$ counts the number of tiles of measure $v_j$ of  type $k$ in the $n$-th iteration of the inflation process. Thus, for the Penrose tessellation with matrix $U^2$ given above, an inflation of the larger tile consists of  1 copy of the smaller and 2 copies of the larger, and an inflation of the smaller tile is the union of 1 copy of each tile type:
$$U^2  
\left(
\ba{c}
1\\
0
\ea
\right) =
 \left(
\ba{cc}
1 & 2 \\
1 & 1
\ea \right)
\left(
\ba{c}
1\\
0
\ea
\right) =
\left(
\ba{c}
1\\
1
\ea
\right) , \quad
U^2  
\left(
\ba{c}
0\\
1
\ea
\right) =
\left(
\ba{c}
2 \\
1
\ea
\right)
$$
and so forth for the powers of $U^d$.
 
\aster
 
Another way to prove theorem \ref{thm: aperiodic silver} depends directly on the positional representation. Let $R'$ denote a fixed   remainder set $R_k$. Since any $x \in R'$ has a representation of the form
$$ 
x = \sum_{k \geq 1} x_{k} \rho^{k}, \quad x_{k} \in \{ 0, 1 \}
$$
each number in the inflated remainder set $\rho^n R'$ can be written $y :=   [  y ] + \{ y \}$
where
$$   [ y ] = \sum_{k=1}^{n}  x_k \rho^{n-k}, \quad \{ y \} = \sum_{k=n+1}^{\infty} x_k \rho^{n-k}  $$

The number $\{ y \}$ is just another element of the remainder set. Regarding $[y]$, if the radix were an integer,  $[y]$ would naturally be called the ``integer part" of $y$. If $\rho$ is not an integer, the set of polynomials in $\rho$ with coefficients in $\{0,1\}$ does not have most of the properties of a ring of integers\footnote{We must adjoin the negative expansions as well.} so we call these numbers  {\it positional representation integers}.

\begin{lemN}
Let $1 < \rho <2$ be a silver number and let $Z[\rho]$ denote the set of positional representation integers, i.e.  the polynomials of finite degree with coefficients in $\{ 0,1 \}$. Then $2 \not \in Z[\rho]$.
\end{lemN}
\vspace{-18pt}
\pf Suppose the contrary. Then 2 is a  polynomial in $\rho$ with coefficients in $\{0 ,1 \}$. Since $\rho >1$, this polynomial cannot have more than 1 term. Hence $2 = \rho^n$ for some $n \in \Z^+$. But this shows $\rho$ is not a silver number. $\Box$ \\

The lemma tells us that the positional representation for $1+1$ is the sum of 1 and a remainder; indeed, the silver equation asserts $2 = 1 \,\cdot \, b_1 \dots b_N$.

Now we can complete the second proof of theorem \ref{thm: aperiodic silver}. If an inflationary tessellation is periodic, then the period lattice is the ring of rational integers. In particular, the sum of any two elements is an element of the lattice.  Since 2 does not belong to the lattice, every tessellation with these tiles is not periodic. Hence the tessellation aperiodic. $\Box$\\

This way of proving aperiodicity extends to higher dimensions and, in principle, also to some tessellations that have just one class of tiles.

\section{2-dimensional positional representations that correspond to aperiodic tilings}

Suppose that $n \in \Z^+$. A consequence of $U v = s v$ is  
\be
U^n v = s^n  v
\label{d-volume}
\ee
We can give this eigenvalue equation a different geometrical interpretation. It asserts that magnification by the factor $s$  of the linear dimensions of a subset of  $\R^{d n}$  of volume  $s v_j$  yields a set whose  volume   is the linear combination $\sum_{k} U_{jk} v_k$ of the $dn$-volumes represented by the $v_k$.

The most accessible cases of this interpretation are $d=1,2$; the most interesting is $d=2$.

Assume that $\rho$ is the golden number $ \phi=\frac{1+\sqrt{5}}{2}  \sim 1.61$. The Frobenius matrix -- the partition matrix --  is 
\be
U = \left(
\ba{cc}
1 & 1 \\
1 & 0
\ea \right)
\label{Ammann-Penrose U}
\ee
Thus
$$
U^2 = \left(
\ba{cc}
2 & 1 \\
1 & 1
\ea \right)
$$
Identify $\R^2$ with the complex field $\C$ and introduce the notation $\langle A \rangle$ for the equivalence class of subsets of  the plane whose area is $A$.

\begin{figure}[h]
\begin{center}
\caption{Radix $\phi$:  Remainder set decompositions based on a non periodic tiling of Ammann.}
\label{Ammann phi n=1}
\includegraphics[width=2.5in]{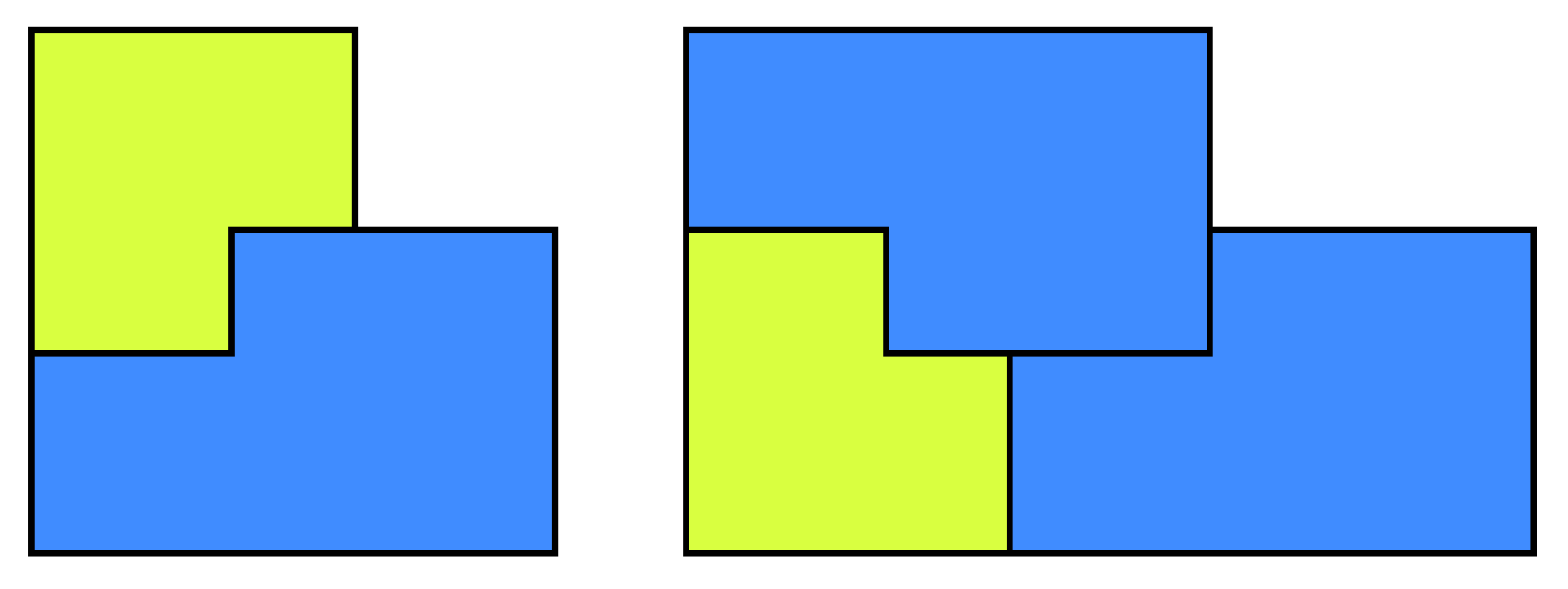}
\end{center}
\end{figure}

\begin{figure}[h]
\begin{center}
\caption{Radix $\phi$: Remainder set decompositions for Penrose remainder sets.}
\label{Penrose 1-digit}
\includegraphics[width=3.5in]{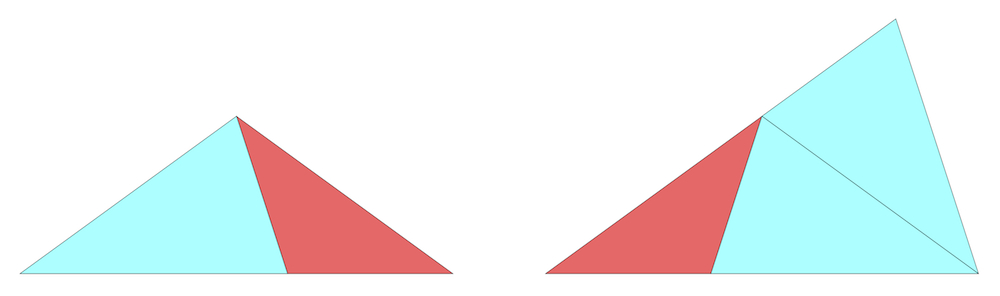}
\end{center}
\end{figure}

The eigenvector of $U$ is $v =  (\langle 1 \rangle, \langle \phi   \rangle )^t$. The relations expressed by $ U^2 v = \phi^2 v$ are precisely those  between the areas of the darts and kites of Penrose's aperiodic tiling of the plane. But not every realization of these relations as areas arises from a conventional tiling because the `tiles' that have the same area need not have the same shape. Figures \ref{Ammann phi n=1} and \ref{Penrose 1-digit} contrast the remainder set decompositions for the aperiodic tilings of Ammann \cite{Ammann1} and Penrose \cite{Penrose1}; both have  radix $\phi$ and the same area condition matrix $U$ (cp. eq(\ref{Ammann-Penrose U})).

\aster
In this context the reader may have observed that the matrix $U$ itself corresponds to the tessellation known as `Ammann's Chair' with radix $\sqrt{\phi}$.\footnote{See fig.\,\ref{Ammann chair 1-digit}  on page \pageref{Ammann chair 1-digit} for the 1-digit decompositions.}   How does this fit in to the general structure? According to what has been  said, $U$ corresponds to a 1-dimensional tessellation. But observe that  $\sqrt{\phi}$ is also a silver number because it is (the largest real) eigenvalue of
\be
U =
\left(
\ba{cccc}
0 & 1 & 0 & 1 \\
1 & 0 & 0 & 0  \\
0 & 1 & 0 & 0\\
0 & 0 & 1 & 0
\ea
\right)
\ee
and
\be
U^2 =
\left(
\ba{cccc}
1 & 0 & 1 & 0\\
0 & 1 & 0 & 1 \\
1 & 0 & 0 & 0 \\
0 & 1 & 0 & 0
\ea
\right)
\label{U:N=4 Amman Chair}
\ee
The characteristic polynomial of the first matrix is the silver polynomial $x^4 = x^2 + 1 $ so of course the roots include the square roots of the golden polynomial.\footnote{Indeed, the $n$-th root of a silver number is also a silver number.} Division by $x^4$ yields the positional representation
$$ 1 = ( 0 \cdot 0\,1\,0\,1 )_{ \sqrt{\phi} } $$
Thus Amman's Chair is a  2-dimensional example  corresponding to the 4-rowed matrix $U^2$ which decouples to describe two tilings, the first and third rows describing the decomposition of  two remainder sets; the third and fourth, the (same) decomposition of the other two. 

\begin{figure}[h]
\begin{center}
\caption{Remainder set $R_{113}$ for the Penrose-like tessellation without complex conjugation. Not edge-to-edge. 6 digits.}
\includegraphics[width=2.5in]{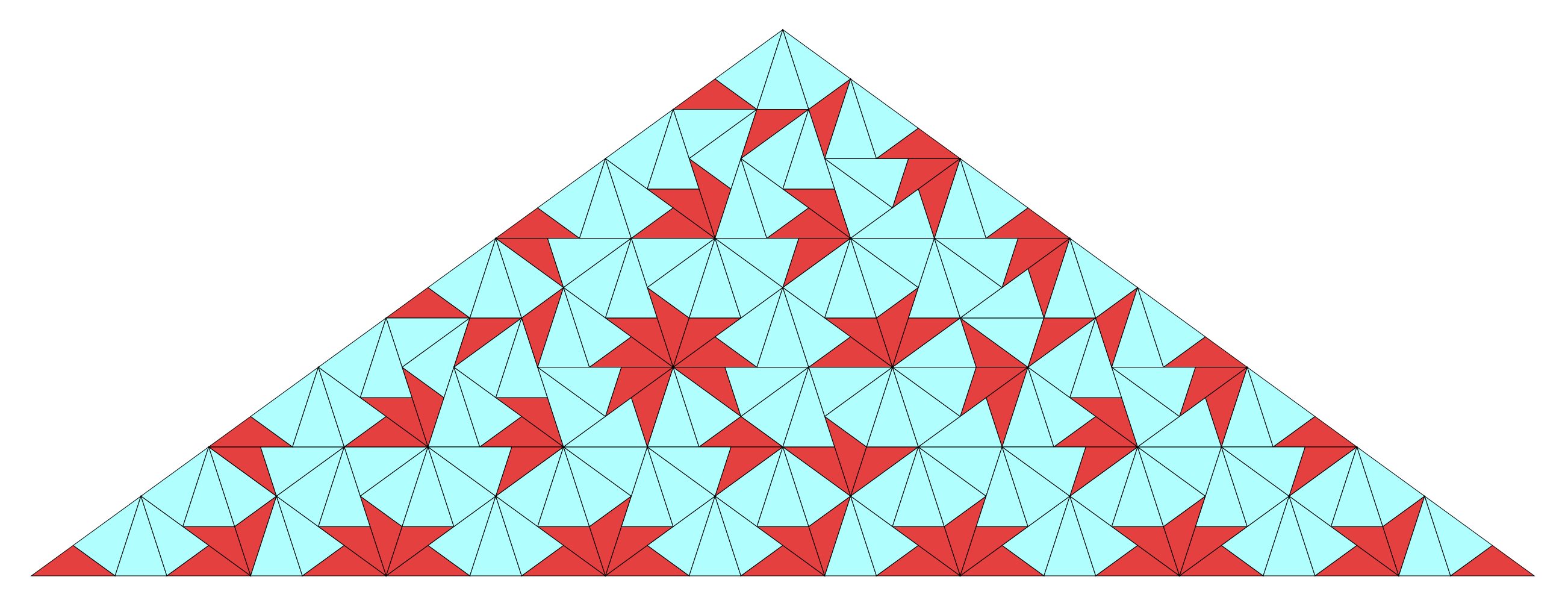}
\label{fig:Penrose R1 no conjugate 6 digit}
\end{center}
\end{figure}

\begin{figure}[h]
\begin{center}
\caption{Remainder set $R_{122}$ for the Penrose-like tessellation without complex conjugation.  Not edge-to-edge. 6 digits.}
\includegraphics[width=2.5in]{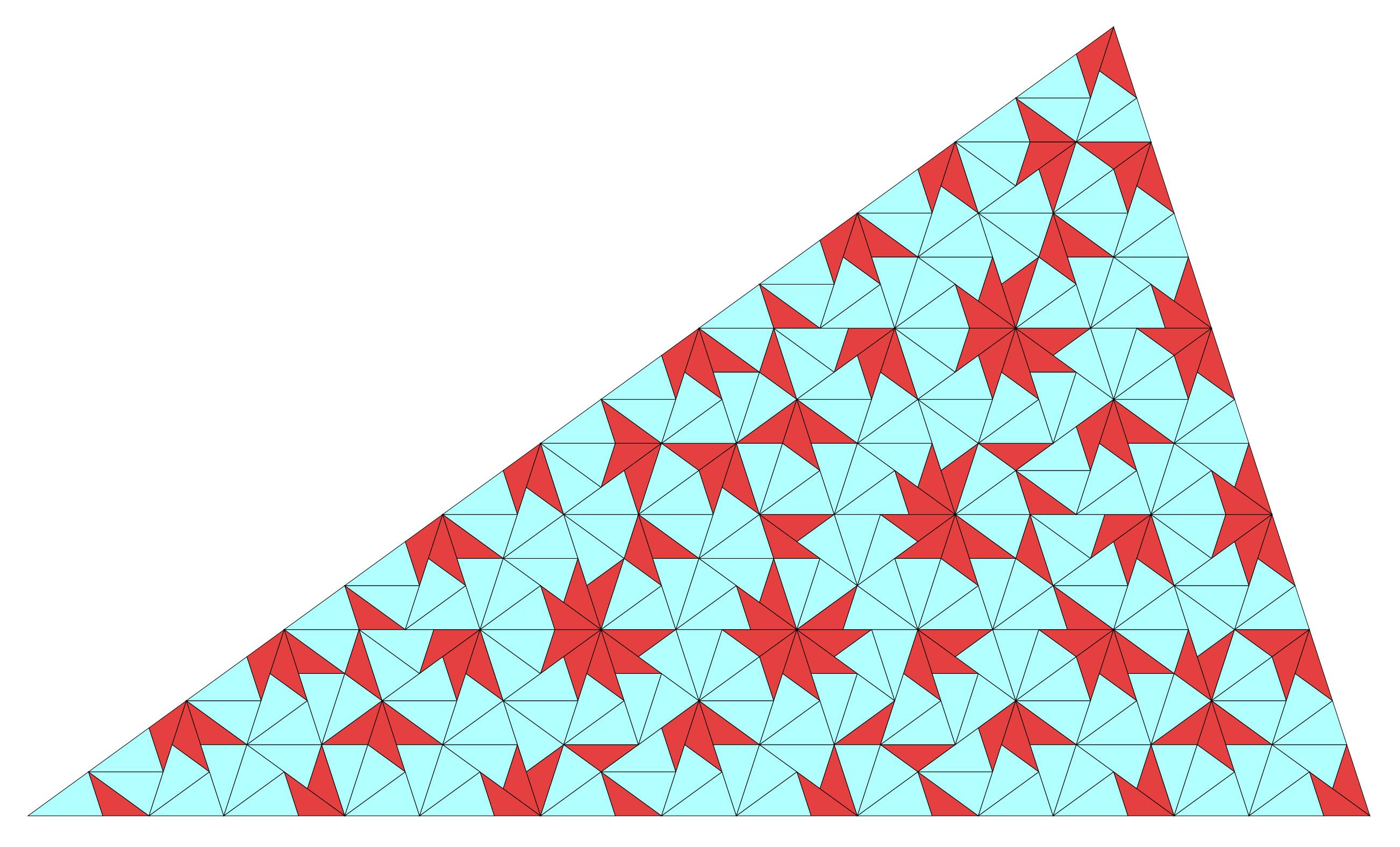}
\label{fig:Penrose R2 no conjugate 6 digit}
\end{center}
\end{figure}

\aster
Consider the group of isometries of $\R^2$ more closely. In the language of complex geometry the orthogonal group $O(2)$ is the extension of the group of rotations by complex conjugation.  Thus one could posit tile equivalence relative to $O(2)$, including conjugation, or merely relative to the proper subgroup of rotations. Assume the latter, and the positional representation with two remainder sets $R_k$ and radix $\rho =\phi$. Abbreviate $u = \exp (i \pi/5)$. An inflationary tessellation is defined by  

\be
\ba{lll}
\rho \, R_{113} &=&( 1 + u^4 R_{113} ) \cup  R_{122}   \\
\rho \, R_{122} &=& ( u/\rho + u^6 R_{113} ) \cup
				 ( 1 + u^3 R_{122} ) \cup
				  ( 1 + u^4 R_{122}  )
\ea
\label{Penrose no conjugate}
\ee
The digits for this positional representation are $\Delta= \left\{ 0, 1, \frac{1} {\rho} e^{i \pi/5} \right\}$. The geometry is  2-dimensional so the number of the digits is the least integer greatest than or equal to $| \rho|^2$. Since $ \phi^2 = \phi +1 \sim 2.61$,   3 digits are expected and that is what we have found. The digits are the coordinates of the vertices of a triangle similar to  $R_{113}$. The partition matrix for this system is the Ammann-Penrose matrix, i.e. the square of  $U$ given by  eq(\ref{Ammann-Penrose U}). This implies that the tiling associated with the recursion eq(\ref{Penrose no conjugate}) describes an aperiodic tiling. Is it the Penrose tiling? The answer is ``No". Although the two remainder sets are the same as the Penrose triangles shown in fig.\,\ref{Penrose 1-digit} and their decompositions also appear to be the same, the tiling described by eq(\ref{Penrose no conjugate}) is constrained by equivalence under the rotation group rather than the full orthogonal group, so fewer  geometrical options are  possible.  In fact, the tessellation associated with eq(\ref{Penrose no conjugate})  is not edge-to-edge. The Penrose edge-to-edge tilings require complex conjugation -- geometrically, reflection -- for their presentation. Figures \ref{fig:Penrose R1 no conjugate 6 digit} and \ref{fig:Penrose R2 no conjugate 6 digit} show   6-digit approximations to the remainder sets without complex conjugation defined by eq (\ref{Ammann-Penrose U}).

\begin{figure}[t]
\begin{center}
\caption{An edge-to-edge Penrose tiling of the triangular remainder set (tile) $R_{113}$ with angles $\{ 1, 1, 3 \} \pi/5$ and radix $\phi$. Tiles that are translations or reflections in the base of the isosceles triangle (complex conjugates) of each other belong to the same type. There are 20 distinct orientation types. 6 digits.}
\label{Penrose R1 20 types 6 digits}
\includegraphics[width=4.5in]{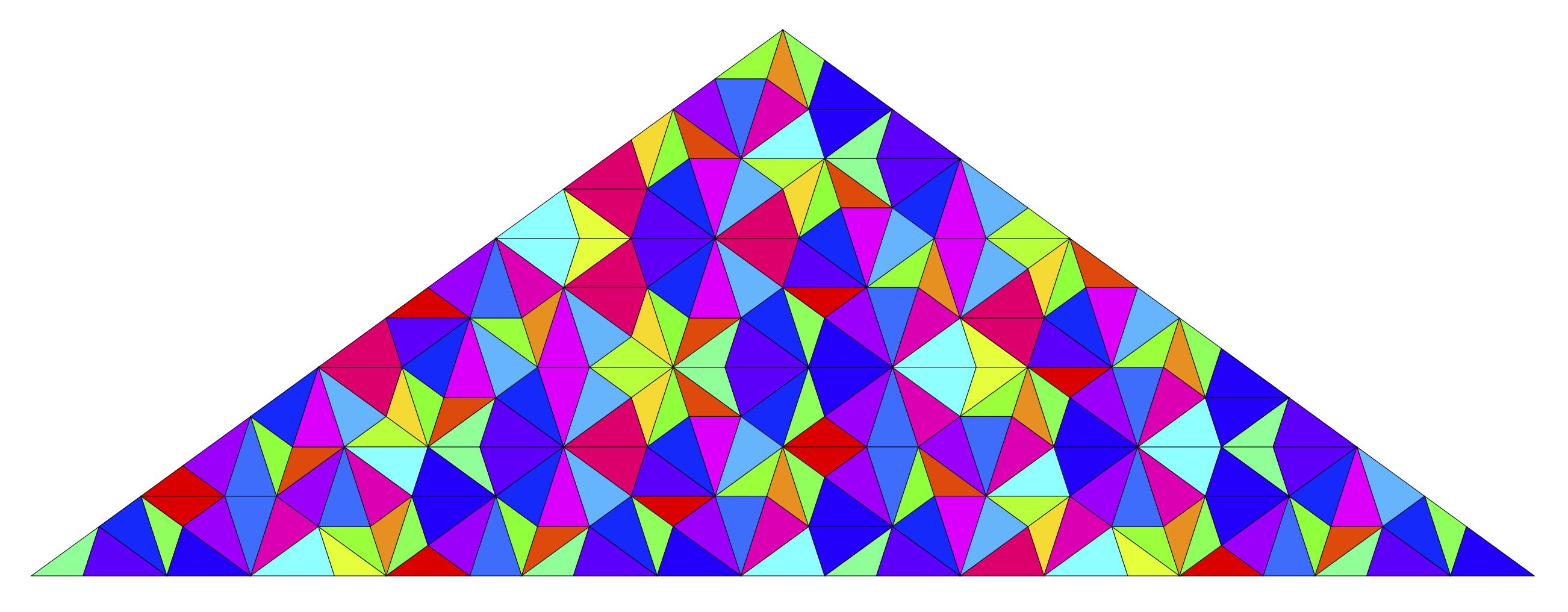}
\end{center}
\end{figure}

\begin{figure}[h]
\begin{center}
\caption{Ammann's Chair: Remainder set decompositions for a positional representation with radix $\sqrt{\phi}$. The small remainder set $R_S$ is shown in red; $R_L$ in blue. $R_L = \sqrt{\phi} \, R_S$.}
\label{Ammann chair 1-digit}
\includegraphics[width=2.0in]{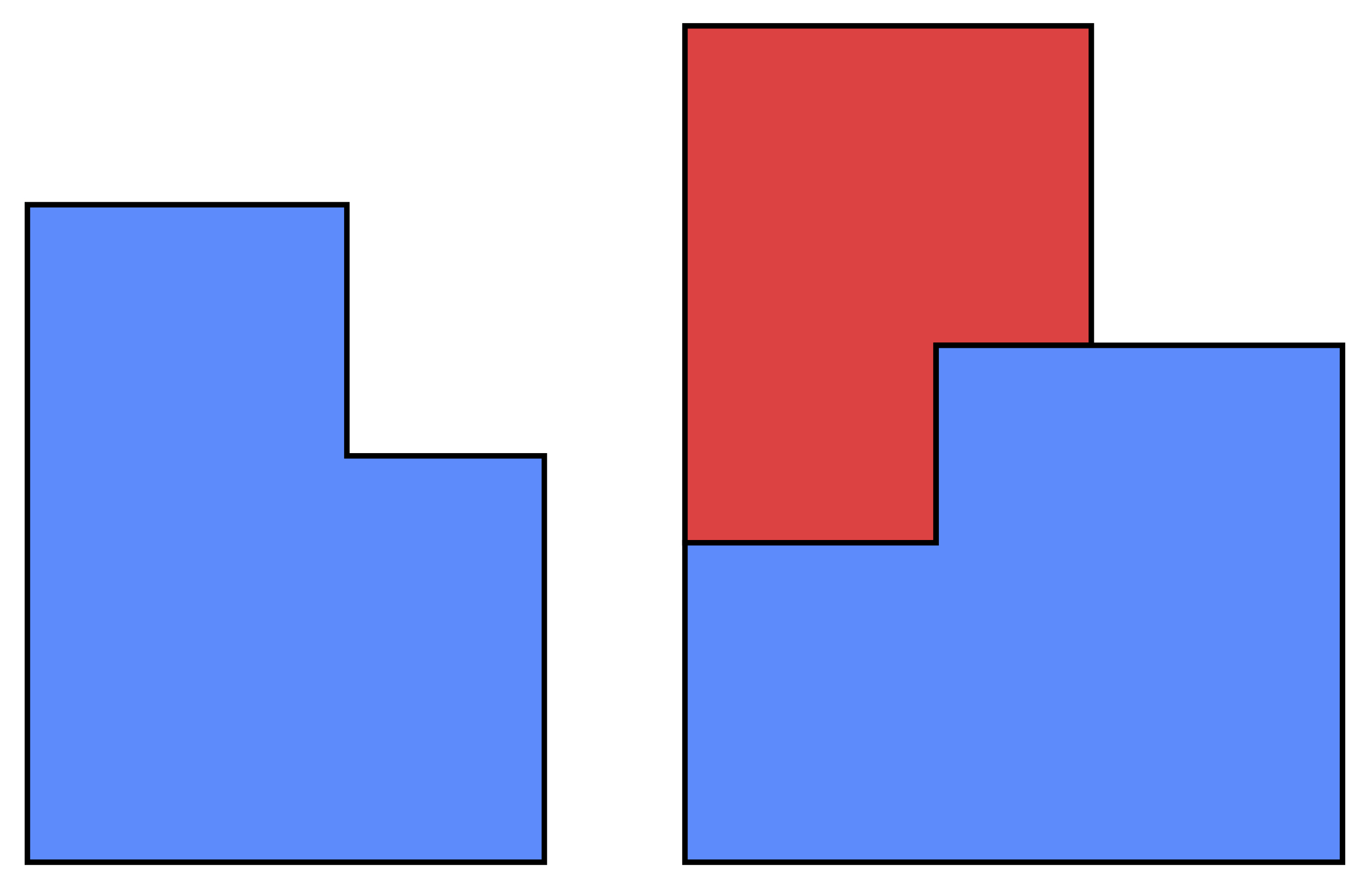}
\end{center}
\end{figure}

Continuing to explore this example, suppose tiles are equivalent if and only if  they are congruent under translation or reflection. Tiles with different orientations  belong to distinct classes, taking into account that orientations and orientations reflected in the real axis are equivalent. Since the rotations are elements of the cyclic group generated by $z \rightarrow e^{i \pi/5} z$, there will be a total of 10 classes for each tile shape: a total of 20. If the tile classes are differently colored, the result is the edge-to-edge pseudo-Penrose tiling shown in figure \ref{Penrose R1 20 types 6 digits} for the triangular  tile with angles $\{ 1, 2,2 \} \pi/5$. This tessellation is aperiodic for the same reason that the usual Penrose tessellation is: the ratios of the number of tiles of different types is irrational.  The partition matrix is a $20 \times 20$ array whose largest eigenvalue is, as expected, $\phi^2$.

Carrying this line of thinking further, suppose that two tiles are equivalent if they are similar hexagons. Then  the `Ammann's Chair' tessellation encountered above, with decompositions shown in figure \ref{Ammann chair 1-digit}, is an aperiodic monotiling of the plane. The `two' tiles are equivalent and hence `the same' since the larger is similar to the smaller (the factor is the radix, $\sqrt{\phi}$; the digits are 0 and 1). Of course, this is not what is normally meant when speaking of a `monotiling',  naive  intuition implicitly expects congruence under the full group of isometries, but the example does emphasize the importance of clearly specifying the equivalence relation when discussing a tiling.

\section{An infinite class of aperiodic tessellations  in $\R^2$}

This section constructs tessellations of $\R^2$ by taking cartesian products of tessellations of $\R$ defined by silver numbers.\footnote{The method can be readily generalized to $\R^d$ but nothing essentially new is gained from the increased complexity of presentation.}

We have already seen that if $U$ is a partition matrix for $\R^1$, then $U^d$ is a partition matrix for $\R^d$. This is how, for instance, the partition matrix for the Penrose and Ammann's Chair tilings arise. However, the partition matrix alone does not determine the geometry of the tiles, as these examples show: the former employs triangles; the latter a non-convex hexagon whose angles are all multiples of $\pi/2$.  

Consider 
\be
1 = \sum_{k=1}^N b_k \rho^{-k}
\label{partition of 1}
\ee
 where $b_k \in \{ 0, 1 \}$ and $b_{N}=1$. The silver number $\rho$ defined by this  equation is the largest real root. The Frobenius companion matrix $U$ for the polynomial is given by eq(\ref{partition matrix}); 
 take it as the partition matrix. The corresponding 1-dimensional geometrical partition of the unit interval consists of the sub-intervals $b_k R_k$ where   $R_k = \{ x : 0 \leq x \leq \rho^{-k} \}, \, 1 \leq k \leq N$. If the $R_k$ are arranged edge-to-edge in any order, their union is essentially disjoint and is an interval of length 1 because of eq(\ref{partition of 1}).\footnote{Some of the $b_k$ may be 0.} Recall that these tessellations of $\R$ are all aperiodic. We will show that their aperiodicity extends the the $d$-dimensional tessellations built up from them. We shall not be concerned with the particular order of arrangement, so suppose it is $R_1, \dots, R_N$ which is the decreasing order of the lengths of the essential subintervals. 
 
Consider the $d$-fold cartesian product of this partition. It provides a partition of the unit hypercube  in $\R^d$ into $n^d$ hyper-rectangles. This partition of the hypercube induces an inflationary tessellation with multiplier $\rho$ for which the hypercube is a remainder set of a $d$-dimensional positional representation with radix $\rho$. We shall prove that it is aperiodic.

Restrict the following discussion to $\R^2$ for simplicity. Define a collection of rectangles by
\be
R_{i,j} = \left\{ (x_1,x_2) : 0 \leq x_{1} \leq  \rho^{-i},  0 \leq x_{2} \leq  \rho^{-j} \right\},  \quad 1 \leq  i, j  \leq N
\label{R_ij}
\ee
The area of $R_{ij}$ is $1/\rho^{i+j}$. Now inflate the largest square $R_{11}$. $\rho R_{11}$ is  a square of side 1. Recalling eq(\ref{partition of 1}), the inflated square can be partitioned into an essentially disjoint union of $\left( \sum b_k \right)^2$ rectangles each of which is a translation and rotation of an $R_{ij}$. Each square $R_{ii}$ occurs once; each $R_{ij}$ with $i<j$ occurs twice. From this information the associated partition matrix $U$, which has order $N(N+1)/2$, is easily constructed. Note that this $U$ is not a Frobenius companion matrix.

\aster

 For the simplest case, $N=2$ so $\rho$ is the golden number. In the previous paragraph it was implicitly assumed that the equivalence relation for tiles is congruence under the full group of isometries.  Now relax this to congruence under translations. Then the number of tile types increases to $N^d$  and the partition matrix entries are 0 or 1. Indeed, the partition matrix for the four tile types is
\be
U =  \left(
\ba{cccc}
1&1&1&1\\
1&0&1&0\\
1&1&0&0\\

1&0&0&0
\ea \right)
\ee

\begin{figure}[h]
\begin{center}
\caption{Radix $\rho = \phi \sim 1.61 $ for the silver polynomial $x^2=x+1$. Six-digit aperiodic decomposition of the  remainder set $R_{11}$ for a positional representation of radix $\rho$.}
\label{dim2 silver(111) 6dig R1decomp}
\includegraphics[width=2in]{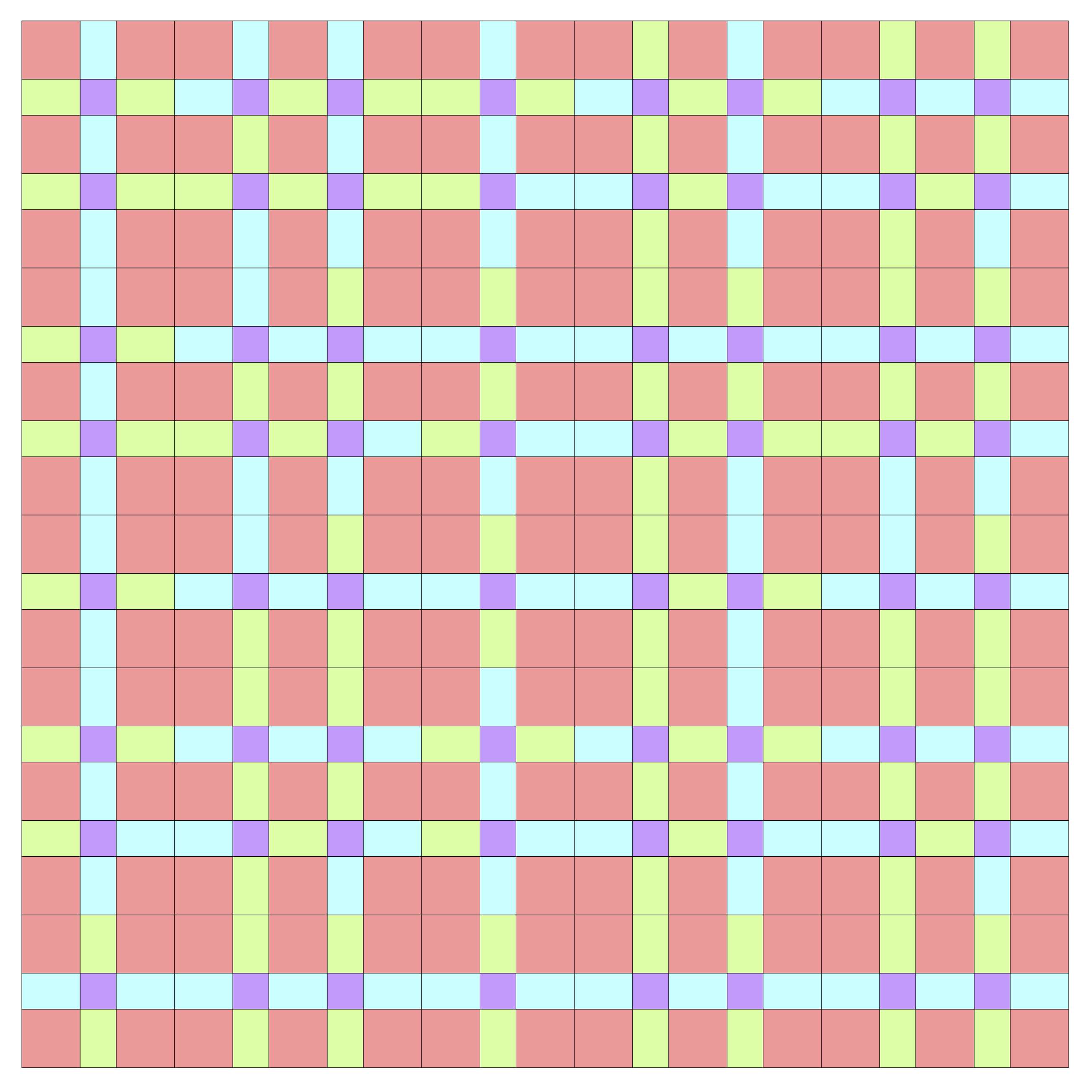}
\end{center}
\end{figure}

The irrational number  $\rho^2$ is the largest eigenvalue; it has eigenvector $(\rho^2, \rho, \rho, 1  )$. The recursion equations are
\be
\ba{lll}
\rho \, R_{11} &=& R_{11}   \cup \left( 1+ R_{12}  \right)  \cup \left(i+ R_{21}   \right) \cup \left( 1+i  + R_{22}   \right) \\
\rho \, R_{12} &=& R_{11} \cup \left( i + R_{21}   \right) \\
\rho \, R_{21} &=& R_{11} \cup \left(1 + R_{12} \right) \\
\rho \, R_{22} &=& R_{11} 
\ea
\ee
The digits are $\Delta = \{ 0,1, i, 1+i \}$. Figure \ref{dim2 silver(111) 6dig R1decomp}  illustrates the 6-digit decomposition of the remainder set $R_{11} = \{z : 0 \leq \Re (z), \Im(z)  \leq 1/\rho \}$.

\aster

Consider $N=3$ and tile equivalence under the full group of isometries. 
The irreducible silver number equation is $x^3=x^2+x+1$ with largest real root $\rho  \sim 1.839$. The partition matrix for the rectangles constructed as the cartesian product of intervals is
\be
U = 
\left(
\ba{cc cc cc}
1 & 2 & 2 & 1& 2 & 1 \\
1 & 1 & 1 & 0 & 0 & 0 \\
0 & 1 & 0 & 1 & 1 & 0 \\
1 & 0 & 0 & 0 & 0 & 0 \\
0 & 1 & 0 & 0 & 0 & 0 \\
0 & 0 & 0 & 1 & 0 & 0 
\ea
\right)
\label{cartesian product U n=3}
\ee
Its largest eigenvalue is $\rho^2$. 
Using these decompositions  to  inflate $R_6$ 4 times produces the tessellation shown in figure \ref{dim2 silver(111) R6 level4}.

\begin{figure}[h]
\begin{center}
\caption{Radix $\rho \sim 1.83 $ for the silver polynomial $x^3=x^2+x+1$. Decomposition of the six remainder sets $R_k$ for a positional representation of radix $\rho$. In order of increasing area, the sets are, from top to bottom and left to right: $R_1$ through $R_6$ (see text). Note that $R_3$ and $R_4$ have equal area. Rotationed and reflected  tiles are equivalent.}
\label{dim2 silver(111) 1dig Rdecomp}
\includegraphics[width=3in]{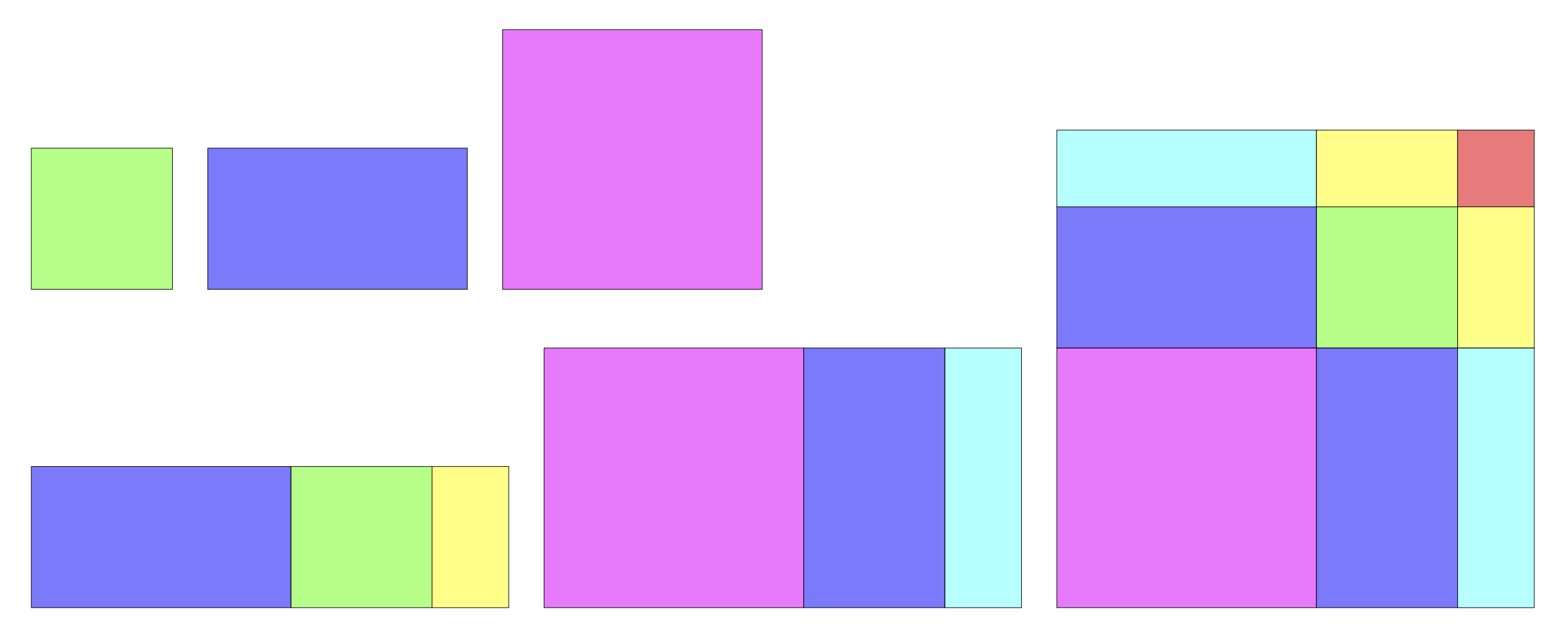}
\end{center}
\end{figure}

\begin{figure}[h]
\begin{center}
\caption{Radix $\rho \sim 1.83 $ for the silver polynomial $x^3=x^2+x+1$. 4-digit tessellation of the largest remainder set $R_6$ for a positional representation of radix $\rho$. See text for details.}
\label{dim2 silver(111) R6 level4}
\includegraphics[width=2.5in]{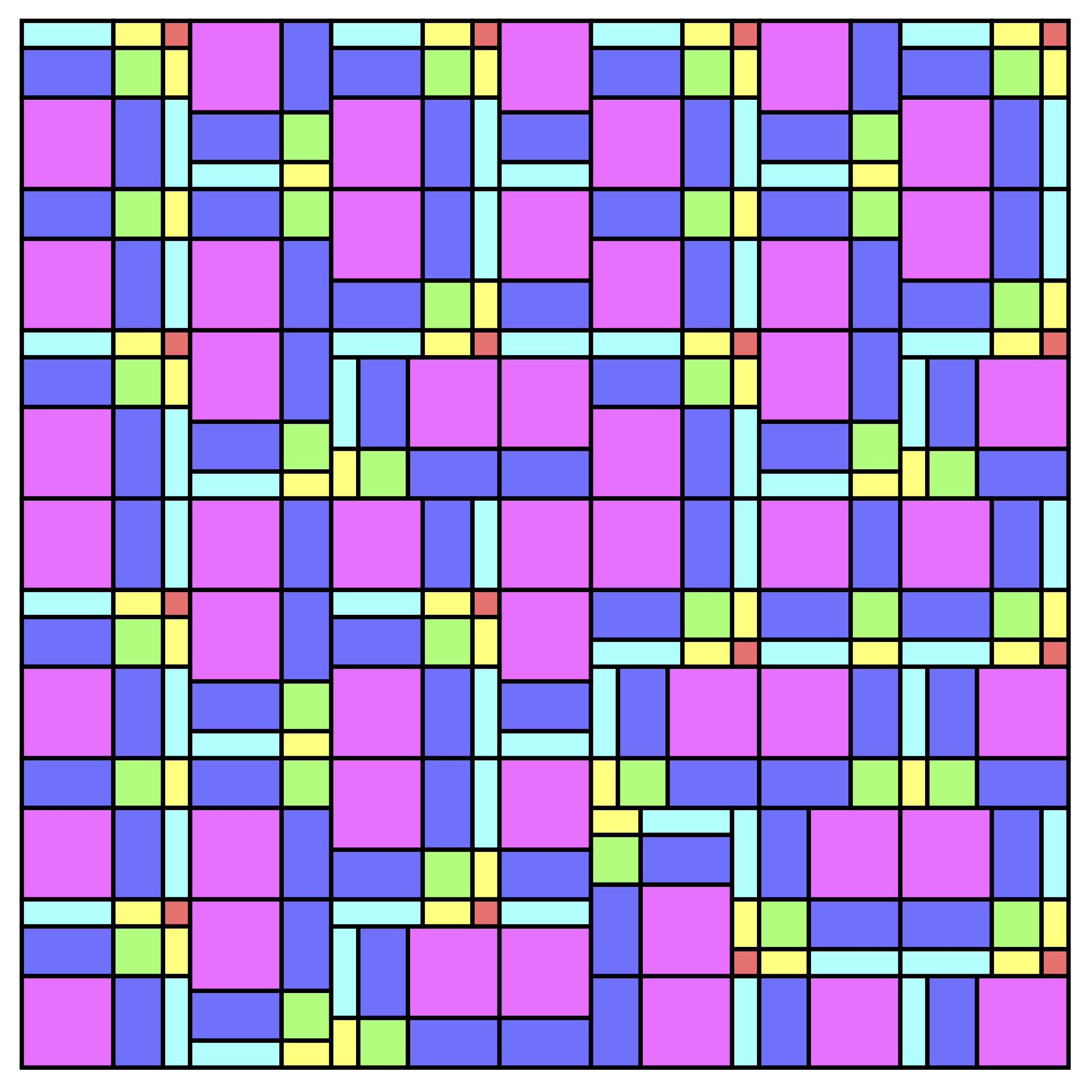}
\end{center}
\end{figure}

\begin{thmN}
The $d$-dimensional tessellations constructed above are  aperiodic.
\end{thmN}
\vspace{-18pt}
\pf The matrix $U$ has $\rho^d$ as largest eigenvalue so Penrose's proof by irrationality applies. $\Box$

\section{Further consideration of tiles equivalent with respect to measure}

At this point we have an extensive inventory of examples and procedures for building 1- and 2-dimensional tessellations. Now it is time to reconsider inflationary tilings for which tiles are said to be of the same type if they have the same measure.

First of all, recall the existence of measure preserving maps $\R^2 \rightarrow \R$.  Identify the unit square $S \subset \R^2$ with the unit interval $I \subset \R$. Divide the square into 4 essentially disjoint congruent squares of side $1/2$ and represent each on one of the 4 essentially disjoint subintervals of length 1/4 that compose $I$; the order doesn't matter. Repetition of this process, subdividing each small square into 4 smaller squares of equal area and representing them on the corresponding  subintervals of equal length, leads in the limit to a mapping of the square onto the line (and, reciprocally, of the line onto the square)  that is well defined except on a set of measure 0. The map  can be extended to a measure-preserving map from $\R^2$ onto $\R$ in the same way. This idea  goes back to the earliest days of Lebesgue measure.

From this description it is evident that such a construction preserves essential disjointness.

\begin{lemN}
Suppose that  $\mu: \R^n \rightarrow \R^{n'}$ is a measure preserving map. Then $A, B \subset \R^{n}$ are essentially disjoint if and only if $\mu(A)$ and $\mu(B)$ are essentially disjoint.
\end{lemN}
\vspace{-18pt}
\pf
Two sets $A$ and $B$ are essentially disjoint if and only if the measure of their intersection is 0, which is equivalent to $m(A \cup B) = m(A) + m(B)$.  If $\mu$ is measure preserving then 
$m(A \cup B) = m'(\mu(A \cup B))$, $m(A)=m'(\mu(A))$, and $m(B)=m'(\mu(B))$ so
$m(A \cup B)  - m(A) - m(B) =0$ if and only if $m'(\mu(A \cup B))  -m'( \mu(A)) - m'( \mu(B)) =0$. $\Box$
\vspace{18pt}

Inflationary tessellations -- positional representations -- provide us with additional explicit structures that establish measure preserving mappings from one dimension to another.  Let us start with an inflationary tessellation of the plane. Denote the collection of remainder sets (or tile class representatives) $R_j$. Let the partition matrix be $U$. Let  $\rho$  be the largest real eigenvalue of $U$. Note that this eigenvalue equation acts on areas, so the eigenvalue is not the multiplier that expands the linear dimensions of a tile. That number is $\sqrt{\rho}$ because the dimension is 2.  

$U$ mediates between the area of an inflated remainder set $m( \rho R_j)$ and the areas of the remainder sets that partition it. The measures $m(R_k)$ are non negative numbers and can be interpreted as lengths. If a remainder set $R_j$ be selected, it can be represented on an integral of length $m(R_j)$, say $I_j:=[0, m(R_j)]$,  and the partition of $R_j$ satisfies
$$ m(R_j) = \frac{1}{\rho} \sum_{k} U_{jk} m(R_k) $$
which describes how many copies of the interval of length $m(R_k)$ are needed to fill out $I_j$. At each level, corresponding to each digit in the representation, the particular order of placement of the sub-intervals within the larger one that contains it, does not matter.  Nor do the translations and rotations and reflections that may occur in the description of the 2-dimensional tessellation matter, for they do not change the measure of the  partitioning tile.\footnote{The measure-preserving maps are functors that carry the structure of inflationary tessellations or positional representations from one space to another.}

Note that the essentially disjoint 2-dimensional inflationary tessellation or positional representation is carried over into a measure-equivalent an essentially disjoint structure by the measure preserving map. This process transfers the structure of the tiling, insofar as it is reflected in the area of each tile, to an interval on the line. We shall call it the {\it projection} of an inflationary  tessellation onto $\R$.

As to the associated radix, note that $U$ refers to the plane tessellation and $\rho$ is the eigenvalue derived from it. After the measure-preserving map is applied and the situation is observed on $\R$, areas have been converted to intervals, the matrix $U$ is still the partition matrix for the nested line segments, and the eigenvalue $\rho$ is -- in this setting -- the multiplier for the inflationary n-dimensional tiling of $\R$. Thus it also is the radix for the positional representation of $\R$ associated to the tiling. 

An interesting example is the Penrose tiling, whose partition matrix is
$$ U = \left(
\ba{cc}
2 & 1\\
1 & 1
\ea
\right)
$$
and largest eigenvalue $\phi^2 \sim 2.61$.  From eq(\ref{Penrose eqs}) on page \pageref{Penrose eqs}, there are just two 1digits' -- 0 and 1. The ratio of the areas $x=m(R_{122})/m(R_{113})$ satisfies $x^2 = x+1$ so $x=\rho$ (because $R_{113} \subset R_{122}$). Choose the unit of length so that  the areas are $(1/\rho, 1/\rho^2)$.

Remembering that the inflated areas increase by the factor $\rho^2$, the two remainder sets can be associated with intervals of lengths $1/\rho$ and $1/\rho^2$, and then partitioned into nested subintervals.  Let the remainder intervals that represent the 2-dimensional remainder triangles be
$$ I_{113} = [0, 1 / \rho], \quad I_{122} = [ 0, 1/\rho^2 ] $$
Then the partition equations can be taken in the form

\be
\ba{lcl}
\rho^2 \, I_{122} &=& \left( \frac{1}{\rho} +  I_{113} \right) \cup   I_{122}  
				 \cup \left( 1 + I_{122} \right) \\
\rho^2 \, I_{113} &=& \left( \frac{1}{\rho} + I_{113}  \right)   \cup I_{122} 
\ea
\label{measure preserv.1}
\ee
with digits $\{ 0, 1, 1/\rho \}$. Recall that $ 1/\rho  = \rho -1$ is an algebraic integer.

\aster

Aperiodicity passes across the dimensional barrier but the issues are slightly different in the two directions. The simpler case is is transfer of a tiling from 2-dimensions to 1. The 1-dimensional `projection'  has the same character as the original.  The reason is simple: a 2-dimensional tiling with equivalence under isometry always projects to a 1-dimensional tiling by line segments with equivalence under isometry -- a tiling in the `usual sense'.

\begin{thmN}
The projection of an aperiodic inflationary tiling is an aperiodic inflationary tiling.
\end{thmN}

In the other direction, the loss of the rotations and reflections as well as arrangement of the positional digits in the plane leaves open  the question of whether there is a 2-dimensional tiling in the sense of equivalence under isometry. If there is, then the character of the original 1-dimensional  tiling `injects' to form an inflationary 2-dimensional tiling.

Thus the 1-dimensional tilings that are projections of the Penrose tiling are also aperiodic, and so forth.

Defining tile type relative to the equivalence relation of equal measure  frees us to investigate a variety of questions about tessellation independent of dimension. The equivalence relation can be constrained later to bring it accord with the conventional view.

\aster

Finally, it may be worth mentioning that there is another way to eliminate the dependence of the tiling equations on rotations and reflections. In a given tessellation, the group generated by the orthogonal transformations that appear in the decomposition equations is finite. Therefore, new tile types and corresponding new remainder sets can be introduced, one for each group element. The original replacement equations imply equations for the new remainder sets but now the replacement equations only contain translations of remainder sets. 

\section{Positional representation integers}

\subsection{$Z[\rho]$ for $\R$}

This section continues the study of 1-dimensional tessellations under the assumptions  that $1< \rho \leq 2$ and  $\Delta =\{ 0,1 \} $. The positional representation integers  are the elements of
\be
Z[\rho] = \left\{  \pm  z :   z =\sum_{k \geq 0}  z_k \rho^k, \quad z_k \in \Delta \right\}
\label{Z(rho)}
\ee
where the sum is finite.  The sum and product of elements of $Z[\rho]$ are defined for only some pairs of elements, but when they are defined, both are commutative and the distributive law $x (y + z) = x \,y + x \,z $ is satisfied for $x, y, z \in Z[\rho]$.\footnote{Such a structure is an example of a {\it ringoid}, but having a name for it is not of much help.}

For  any $z \in Z[\rho]$ and $n \in \Z$, 
$$ 0 + z = z, \quad 1\, z =z,  \quad \rho^n Z[\rho] \subset Z[\rho] $$
so $0$ is an additive unit and $1$ a multiplicative unit, and inflation of $Z[\rho]$ by the multiplier $\rho$ is a subset of $Z[\rho]$. Each element has an additive inverse.  If $\{ m_j \} \cap \{ n_k \} = \emptyset$ then $\sum  \rho^{m_k} + \sum \rho^{ n_k}$  is always defined whatever $\rho$. Denote the number of elements in a finite set $S$ by $\#(S)$. If $\#\left( \{ m_j + n_k \} \right) = \#\left( \{ m_j  \} \right) *  \#\left( \{  n_k \} \right)$  (i.e., the exponents of the termwise products are all different) then $ \left( \sum  \rho^{m_k} \right) \left(  \sum \rho^{ n_k} \right)$ is defined.

In general, $2 = 1 + 1  \not \in Z[\rho]$. Of course 2 always has a positional representation but a remainder may appear. 

The relationship between the sets of integers $Z[\rho]$ and lattices will be considered next. Three examples in $\R$ suggest how different these sets of `integers' can be depending on the radix. In all cases the set of digits is $\Delta$. A lattice $\Lambda \subset \R$ has a single generator so $\Lambda = \lambda \Z$ for some $\lambda >0$ in $\R$. Unless otherwise stated, assume that $1 < \rho <2$ and $\Delta = \{ 0, 1 \}$.

\begin{propN}
If $\rho=2$ then $Z[\rho]  = \Z$ is a lattice and a ring. 
\end{propN}

This positional representation is just the signed standard binary representation for integers.

\begin{propN}
If $1 < \rho < 2$ then $Z[\rho]$ is not contained in a lattice.
\end{propN}
\vspace{-18pt}
\pf Suppose $Z[\rho]  \subset \Lambda$. Then $1 \in Z[\rho] $ implies $1= n \lambda$ with $n \in \Z^+$. Hence $\lambda = 1/n$. Similarly,   $\rho^j  = k_j/n$ so $m^j = k_j n^{j-1}$ for $j \geq 1$ so $n^{j-1}$ divides $m^j$. Let the prime factorizations be
$$ m = \prod_l p_l^{d_l}, \quad n= \prod_l p_l^{e_l} $$
The divisibility condition is
$$  k_j = \prod_l p_l^{j d_l - (j-1)e_l}$$
so $ d_l \geq  (1-1/j)e_l $ for $j>1$. Taking $j>>1$ and noting that $d_j$ and $e_j$ are integers, it follows  that $d_j \geq e_j$.  At least one of these inequalities must be strict; otherwise $m=n$, a contradiction. But $d_j>e_j$ implies $\rho = m/n \geq p_j \geq 2$, a contradiction.
$\Box$
\vspace{18pt}

Now consider whether a lattice can be a subset of $Z[\rho]$. Here is an interesting example.
 
\begin{propN}
If $\rho=\sqrt{2}$ then $Z[\rho]$ contains $ \Z$ as a proper subset.
\end{propN}
\vspace{-18pt}
\pf Every integer $n \geq 0$ in $Z[\rho]$  has a representation of the form
\bean
 n &=&  \sum_{k \geq 0} n_k 2^{k/2} \\
 &=&   
\sum_{k \geq 1}  n_{2k} 2^k  + \sqrt{2} \sum_{k \geq 0}  n_{2k+1} 2^k\\
&=& m1 + m_2  \sqrt{2}, \quad m_1, m_2 \, \mbox{non negative rational  integers}
\eean
In particular, $\Z \subset Z[\rho]$ and  $ \sqrt{2}\,  \Z \subset Z[\rho] $.  Both inclusions are proper because $\sqrt{2} \not \in \Z$ and $1 \not \in \sqrt{2} \,\Z$. Moreover,  $Z[\rho]$  is dense in $\R$, hence not a subset of a lattice. $\Box$
\vspace{18pt}

This example is an instance of an interesting class of positional representations. The remainder set is $R=[ 0,1+\sqrt{2} ]$.  Each $x \in R$ can be written in the form $x= x_1 + \sqrt{2} \, x_2$ for $x_1, x_2 \in [0,1]$. Each $x \in (1, \sqrt{2} )$ has uncountably many different representations because $x_1 \in (0,1)$ can be chosen at will whence  $x_2 = (x-x_1)/\sqrt{2} \in (0,1)$. 

\begin{propN}
If $\Lambda = \lambda \Z  \subset Z[\rho]$ then $\lambda$ and $\rho$ are of the same type, i.e.  algebraic or transcendental. If $\rho$ is algebraic, then both belong to the number field generated by $\rho$ and are algebraic integers.
\end{propN}
\vspace{-18pt}
\pf For each  $k \in \Z^+$ there is a $Q_k \in Z[\rho]$ such that $ k \lambda = Q_k $. Each $Q_k$ is a polynomial in $\rho$ whose coefficients are 0 or 1. In particular, $\lambda = Q_1$. It follows that $\lambda$ and $\rho$ are of the same type, i.e. if one is, respectively,    algebraic or transcendental, so is the other. 

  From $k = \frac{k \lambda}{\lambda} =Q_k/Q_1$ it follows that  $Q_k -k Q_1=0$ and the polynomial has leading coefficient 1 (Note that if $k>l$ the degree of $Q_k$ is greater than the degree of $Q_l$).  Hence $\rho$ is an algebraic integer. Therefore $k \lambda =Q_k$ lies in the field generated by $\rho$. $\Box$
\vspace{18pt}
 
If a plane tessellation is periodic and inflationary, then   theorem \ref{thm 1} implies that $\rho$ is a quadratic imaginary integer and $\rho \, O (\Lambda) \subset \Lambda$ where $O$ is a suitable orthogonal linear transformation. Continuing to restrict our attention to multipliers  for which $1 < \rho^d \leq 2$, which means here that $1 < \rho^2 \leq 2$, the suitable  quadratic imaginary fields are given by $\rho \,  O= i \sqrt{2}$, $\rho= 1+i $, and $\rho= \frac{1+i \sqrt{7}}{2}$ and their multiplies by units of the field $\Q(\rho)$ each generates. In each case $Z[\rho]$ is the associated ring of algebraic integers, which is a lattice.

\section{Monotile aperiodic tessellations}

Joan M. Taylor recently devised a monotile -- a tile that can tessellate the plane only aperiodically \cite{Taylor 2010}, \cite{Socolar+Taylor 2010}, \cite{Socolar+Taylor 2011.2}.
 
\begin{figure}[t]
\begin{center}
\caption{Radix = 2. Taylor aperiodic monotile. The line segments  have  zero width and are shown only for illustrative purposes. The tessellation is not inflationary.}
\label{Taylor monotile}
\includegraphics[width=2in]{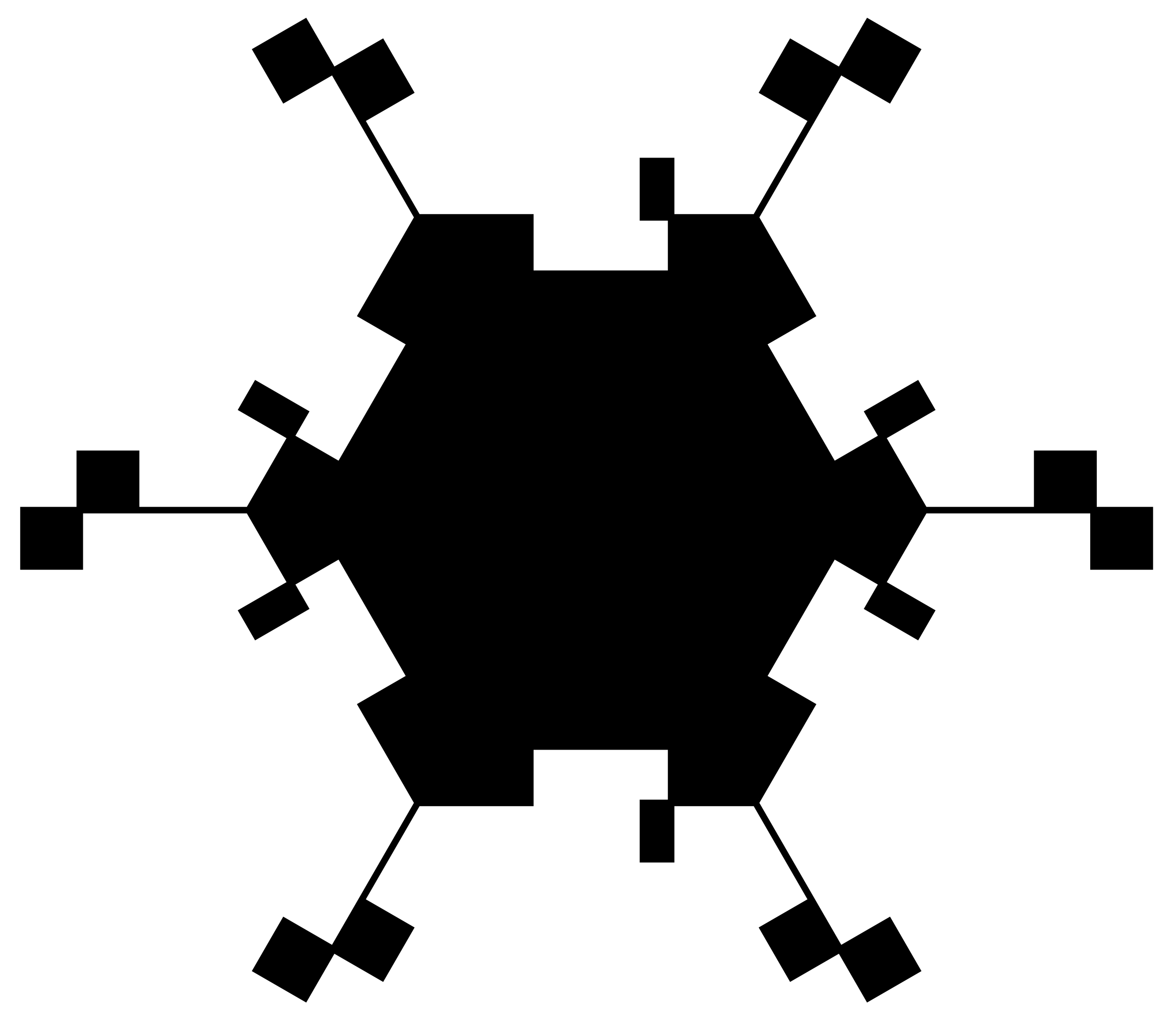}
\end{center}
\end{figure}

There are various ways to think about this tiling. From a conventional  perspective, it is the standard hexagonal tiling of the plane but the tiles are `decorated' and replacement rules are promulgated that insure only aperiodic tilings arise. There is only one type of decoration. This description -- which is standard in the field -- may give the impression that not so much has been accomplished; some might be better satisfied with a specific geometric shape that only allowed aperiodic tilings.  Socolar and Taylor have shown that such a tile exists: the decorations are equivalent to a hexagon whose boundary has been suitably deformed. The deformation introduced by Socolar and Taylor is shown in figure \ref{Taylor monotile}. Note that this tile is not self-similar.\footnote{See the next subsection for a self-similar monotile.}

Taylor's monotile, shown in figure \ref{Taylor monotile},  is a distorted version of a regular hexagon.\footnote{This image agrees with  the black tile in Figure 6(a) of \cite{Socolar+Taylor 2011.2} but differs from what purports to be the same tile in figure 3 of \cite{Socolar+Taylor 2010}.}  The area of the monotile is the same as the area of the hexagon from which it was constructed.

\begin{figure}[t]
\begin{center}
\caption{The Taylor monotile yields an essentially disjoint covering of the plane.}
\label{Taylor monotile 1dig}
\includegraphics[width=2in]{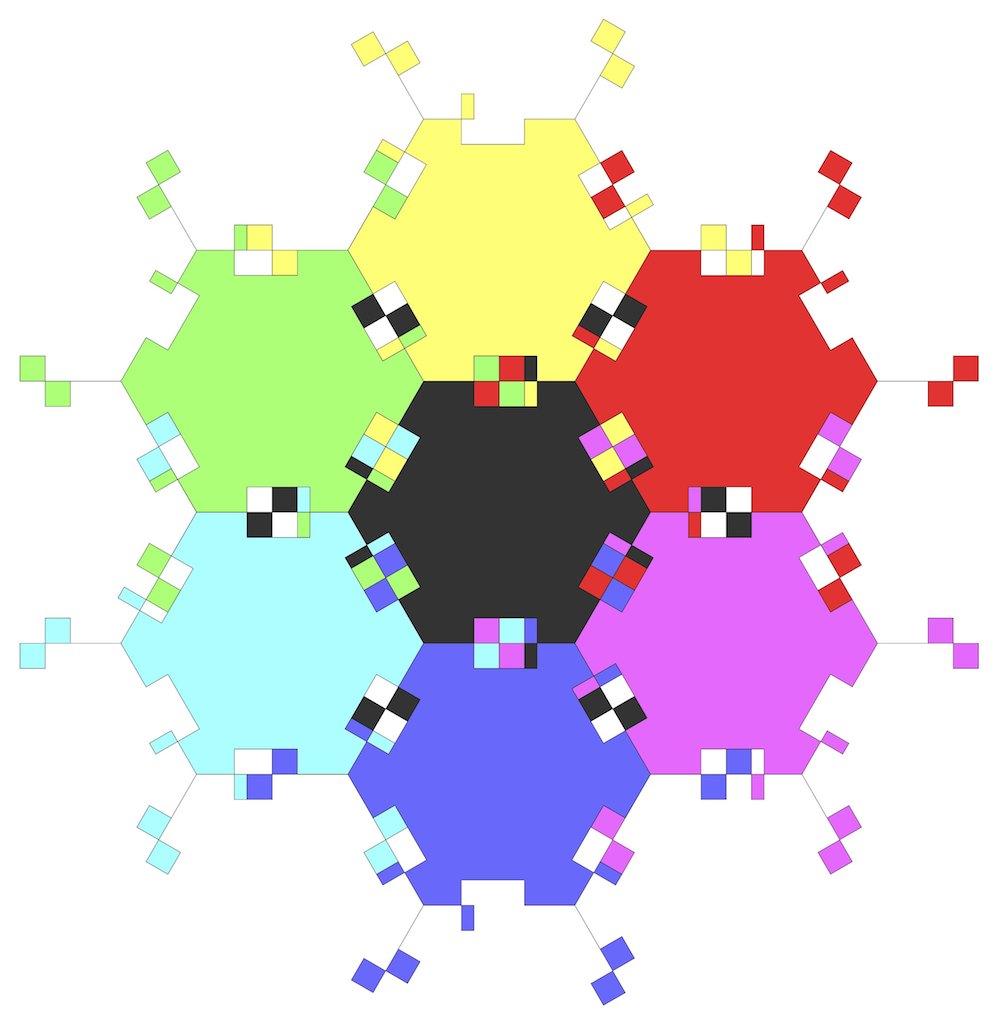}
\end{center}
\end{figure}

Figure \ref{Taylor monotile 1dig} shows 7 interlocking copies the Taylor monotile so the monotile is not self-similar. Socolar and Taylor prove  that this process can be continued to an essentially disjoint covering the plane and that the covering is aperiodic. The figure employs color to make it easier to see how the 7 copies of the monotile interlock, but in this figure color does not play an intrinsic role -- it does not carry any  information about the tessellation additional to tile's shape.  The union of the 7 tiles also has a generally hexagonal shape, and it follows that if the process is repeated so that there are $7^n$ copies of the monotile after which  the scale of the figure is divided by $2^n$, then the limit will be a regular hexagon whose side can be taken to have length 1. 
This presentation is also equivalent to employing congruent hexagonal tiles that differ only in color. There will be 4 types, say red, green, teal and magenta. 

The tessellation is not inflationary and therefore it is not an example of the tilings considered in this paper.  However, a second tiling due to Taylor, whose discovery chronologically preceded the tiling just described, is inflationary. We shall study this  tiling in more detail.

\aster

The method of proving a tiling is aperiodic by showing that the limit of the ratio of the number of tiles of different types that occurs in a sequence of regions of increasing and unbounded measure is irrational cannot be used when there is only one type of tile.  

An approach which would apply to monotilings might be to decorate each tile type with a line segment so that the partition rules for tile placement insure that the line segments become continuous curves without endpoints. As a tile is inflated and the rules are applied, increasingly large regions of the plane will be covered. The condition that a curve constructed this way not have endpoints  means it starts and ends at infinity, or is a closed curve in a compact region of the plane. If there were a nested collection of closed curves of increasing length, then the tessellation could not be periodic. 
This is the principle that Taylor \cite{Taylor 2010} and Socolar and Taylor \cite{Socolar+Taylor 2011} use to prove the existence of an aperiodic monotile, which will be discussed below. But it is well to note that it could also apply when there is more than one type of tile. 

\begin{figure}[h]
\begin{center}
\caption{Decorated square tessellation showing nested sequence of concentric squares of increasing size. See text for details.}
\label{square lines}
\includegraphics[width=1.5in]{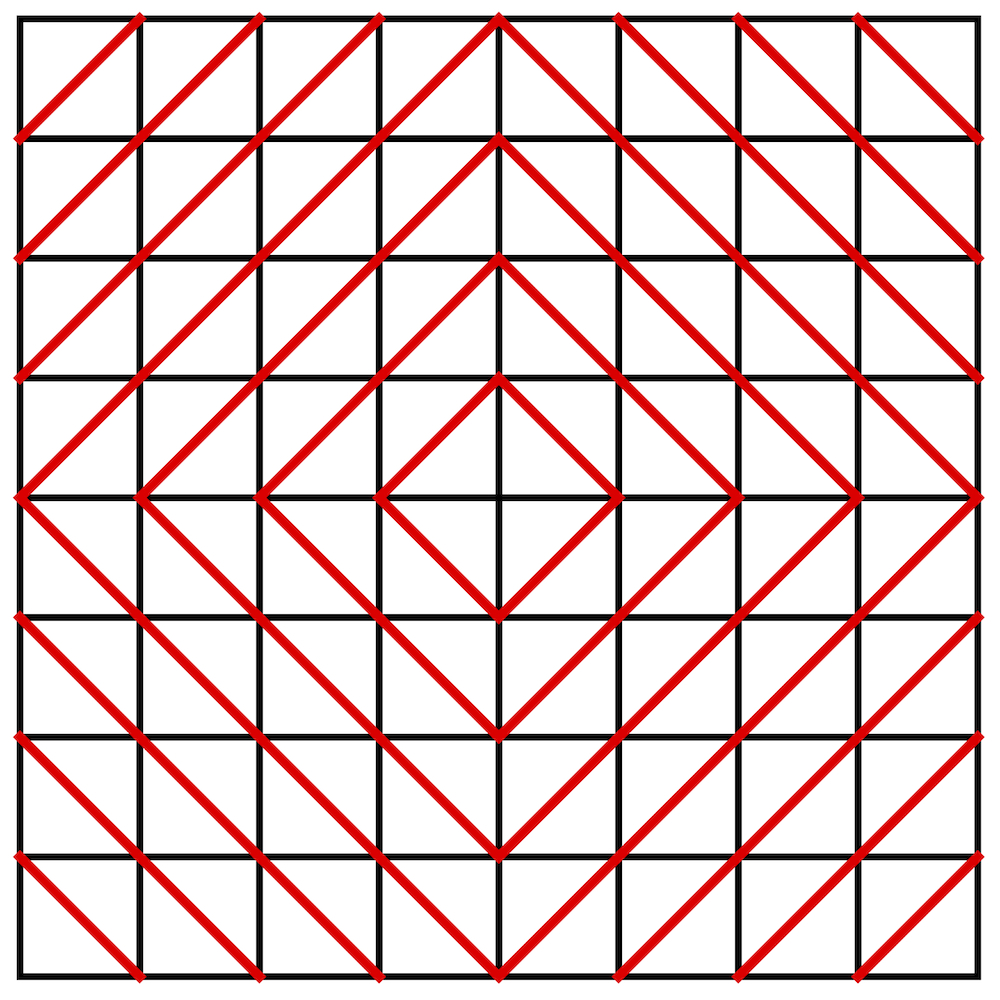}
\end{center}
\end{figure}

Although the result is valid, the principle as stated is not correct. Consider, for example, the inflationary tiling of the plane by translated squares defined by eq(\ref{squares})  on page \pageref{squares}.  After inflation, the recursion fills the first quadrant. Extend the covering to the plane by rotating the figure by multiples of $\pi/2$.  There is only one type of tile and there are no rotations. Decorate the square with a diagonal line. The decoration of the basic tile propagates in the tessellation, shown in fig.\,\ref{square lines}, as a sequence of concentric nested squares of unbounded size. These curves are continuous, closed and have finite length. The emergent red pattern is not periodic, but the underlying geometric tiling is periodic. 

The kernel idea that nested curves of increasing length imply non periodicity is basically correct but it requires some modification. If the underlying tiling were periodic with period lattice $\Lambda$, then decorations of the tessellation would project onto the quotient torus $\C/\Lambda$. It is the nature of these curves that decides the question of periodicity. If the tessellation is periodic, then the images of the decorations on the torus will be a (system of) closed curves of finite length. Otherwise, the tessellation cannot be periodic.

Returning to the decorated square, the projected curves collapse to a simple loop on the torus.

Returning to the Penrose tessellation, recall that it has an inflationary construction based on two remainder sets: triangles $R_{jkk}$ where a subscript, say $j$, denotes the angle  $j \pi/5$. 

\begin{figure}[h]
\begin{center}
\caption{Penrose pinwheel based on remainder set $R_{122}$ with line segments, showing nested continuous curves of increasing length. 5 digits. The image is generated by the positional representation eq(\ref{Penrose eqs}). See text  for details.}
\label{Penrose R122 pinwhl lines}
\includegraphics[width=3in]{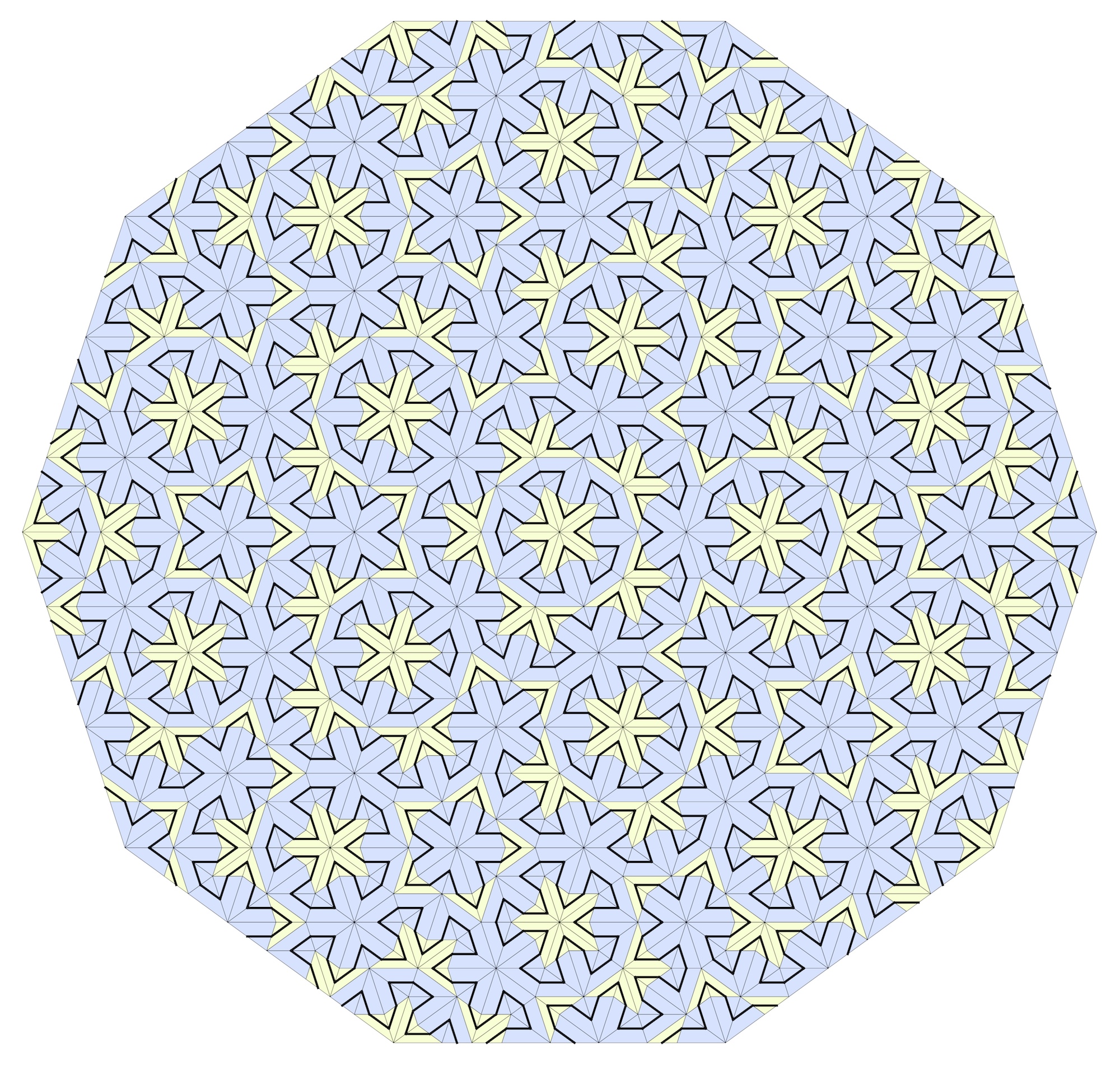}
\end{center}
\end{figure}

The remainder sets $R_{113}$ and $R_{122}$ are isosceles triangles. Mark $R_{113}$ with a line segment joining the midpoints of the isosceles edges, and $R_{122}$ with a line segment joining the midpoint of the base to the opposite (acute) vertex. Figure \ref{Penrose R122 pinwhl lines} shows a Penrose 5-digit pinwheel tiling based on $R_{122}$ where only the line segments are indicated. They form continuous curves. It is easily seen (but it is not quite so easy to prove) that there are nested sequences of closed curves that meander about the perimeter of a sequence of regular pentagons of increasing size centered on the axis of the pinwheel. It seems intuitively clear that the curves cannot project onto closed curves of finite length on any quotient torus, but it may not be easy  to prove this.

\aster

The high degree of interest in the recently discovered aperiodic monotile -- a single tile class that yields only aperiodic tilings -- may make it worthwhile to to provide more details about its connection with positional representation.  Taylor's initial discovery involved 14 congruent trapezoids  with decorations that determined their placement. This was later supplemented by the equivalent (disconnected) monotile described above and illustrated in figure\,\ref{Taylor monotile}.

Fourteen congruent trapezoids are suitably decorated so that appropriate tiling rules pair them to form hexagons. A hexagon inflated by the radix $\rho=2$ can be tiled by four  decorated congruent hexagon tiles and 6 congruent trapezoids. Repeated inflations cover an ever increasing internal region with hexagons and a region contiguous to the boundary of the inflated hexagons with trapezoids. We will derive the equations that correspond to Taylor's trapezoidal construction. 

Reference \cite{HLR:mono tile}  gave a system of coupled equations for the 14 decorated Taylor trapezoidal remainder sets as well as a single equation that encapsulates all the information. Next we present an equivalent but slightly simpler version of that  equation which describes the aperiodic  Taylor tiling.\footnote{The system described below orients the trapezoid vertically and has simpler digits.}

Put $\rho=2$ and $\omega=e^{ \pi i/3}$. Denote complex conjugation by an overbar.  

\begin{thmN}
 The recursion eq(\ref{HLR mono trap recursion})  produces the Taylor-Socolar aperiodic  monotiling. The digits are $\{ 0, 1 , \omega, \overline{\omega} \}$. 
\be
\rho \, R = R \cup 
	   \left( \omega  + \overline{\omega}^2 R   \right) \cup
	 \left( \overline{\omega}  + \omega^2 \, \overline{R}   \right) \cup
	  \left( 1  + \omega^3  \,  \overline{R}   \right) 
\label{HLR mono trap recursion}
\ee
\end{thmN}

\begin{figure}[t]
\begin{center}
\caption{A portion of the trapezoid-based monotiling produced by eq(\ref{HLR mono trap recursion}).  The field of trapezoids is shown in the background for reference. Hexagons, related to the Socolar-Taylor hexagonal monotile,  are formed by pairs of trapezoids. 5 iterations.}
\label{HLR mono trap 5 dig}
\vspace{4pt}
\includegraphics[width=3in]{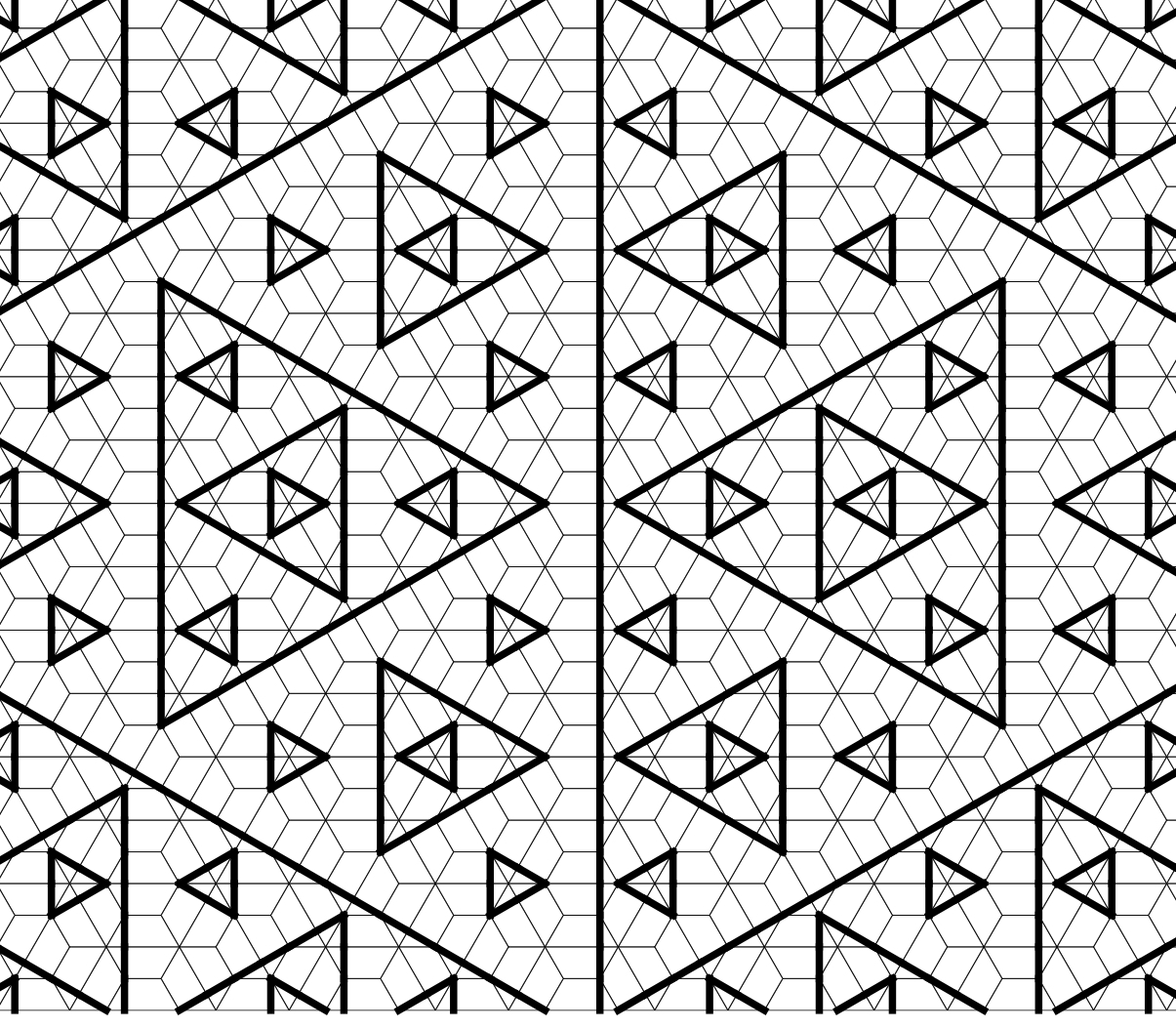}
\end{center}
\end{figure}

Figure \ref{HLR mono trap 5 dig} illustrates Taylor's aperiodicity proof. The figure was constructed from eq(\ref{HLR mono trap recursion}) as follows: Draw a line perpendicular to the parallel edges of the trapezoid from a vertex to the base. The resulting figures consists of an array of triangles. The idea of the proof of aperiodicity  is that the triangles have increasing size as the tessellation is inflated to cover more of the plane, and that arbitrarily large triangles are incompatible with periodicity. 

The alternative argument  based on the fact that a fundamental domain for a periodic array can be translated by lattice elements without changing the lattice can in principle be applied to monotile tessellations.  If the array of positional integers, i.e. the set of polynomials whose coefficients are digits drawn from $\Delta$ evaluated at the radix $\rho$, is not a lattice then the tessellation by remainder sets cannot be periodic. This method does apply to the Penrose aperiodic tiling and the product space tilings with silver number radix, but it does not apply to the Taylor trapezoid monotile for in that case, the set of positional integers is the lattice $\{ z = m +\omega n \}$ where $\omega = \exp (i \pi/3)$.

\end{document}